\pdfoutput=1
\RequirePackage{ifpdf}
\ifpdf 
\documentclass[pdftex]{sigma}
\else
\documentclass{sigma}
\fi

\numberwithin{equation}{section}
\newtheorem{Theorem}{Theorem}[section]
\newtheorem{Corollary}[Theorem]{Corollary}
\newtheorem{Lemma}[Theorem]{Lemma}
\newtheorem{Proposition}[Theorem]{Proposition}
{ \theoremstyle{definition}
\newtheorem{Definition}[Theorem]{Definition}
\newtheorem{Example}[Theorem]{Example}
}

\DeclareMathOperator{\rank}{rank}
\DeclareMathOperator{\Span}{Span}

\begin{document}

\newcommand{\arXivNumber}{1312.4581}

\allowdisplaybreaks

\renewcommand{\PaperNumber}{031}

\FirstPageHeading

\ShortArticleName{Invariants and Inf\/initesimal Transformations for Contact Sub-Lorentzian Structures}

\ArticleName{Invariants and Inf\/initesimal Transformations\\
for Contact Sub-Lorentzian Structures\\
on 3-Dimensional Manifolds}

\Author{Marek GROCHOWSKI~$^{\dag\ddag}$ and Ben WARHURST~$^\S$}

\AuthorNameForHeading{M.~Grochowski and B.~Warhurst}

\Address{$^\dag$~Faculty of Mathematics and Natural Sciences, Cardinal Stefan Wyszy\'{n}ski University,\\
\hphantom{$^\dag$}~ul.
Dewajtis 5, 01-815 Waszawa, Poland} \EmailD{\href{mailto:m.grochowski@uksw.edu.pl}{m.grochowski@uksw.edu.pl}}

\Address{$^\ddag$~Institute of Mathematics, Polish Academy of Sciences,\\
\hphantom{$^\ddag$}~ul.~\'{S}niadeckich 8, 00-950 Warszawa, Poland}
\EmailD{\href{mailto:m.grochowski@impan.pl}{m.grochowski@impan.pl}}

\Address{$^\S$~Institute of Mathematics, The Faculty of Mathematics, Informatics and Mechanics,\\
\hphantom{$^\S$}~University of Warsaw, Banacha 2, 02-097 Warszawa, Poland}
\EmailD{\href{mailto:benwarhurst@mimuw.edu.pl}{benwarhurst@mimuw.edu.pl}}

\ArticleDates{Received October 10, 2014, in f\/inal form March 30, 2015; Published online April 17, 2015}

\Abstract{In this article we develop some elementary aspects of a~theory of symmetry in sub-Lorentzian geometry.
First of all we construct invariants characterizing isometric classes of sub-Lorentzian contact $3$~manifolds.
Next we characterize vector f\/ields which generate isometric and conformal symmetries in general sub-Lorentzian manifolds.
We then focus attention back to the case where the underlying manifold is a~contact $3$~manifold and more specif\/ically
when the manifold is also a~Lie group and the structure is left-invariant.}


\Classification{53B30; 53A55; 34C14}

\section{Introduction}

\subsection{Basic notions and motivation}

Sub-Lorentzian geometry is a~relatively new subject although it does fall within the scope of broader perspectives on geometry.
For instance the work of Berestovskii and Gichev~\cite{BerGic} on metrized semigroups is perhaps the broadest
perspective one could take on the subject or alternatively a~more closely related generalisation is the notion of
para-CR-geometry see~\cite{HN}.

The aim in this work is to look specif\/ically at sub-Lorentzian geometry and so in this section we only present the
notions that are required for the formulation of the main results of the paper.
For more details and facts concerning sub-Lorentzian geometry, the reader is referred to~\cite{g6} and the references
therein (see also~\cite{Grong, markin}).
The structures which support the sub-Lorentzian structures in this paper are identical to those which support the
analogous structures in the sub-Riemannian setting and so the reader is also referred to the paper of Agrachev and
Barilari~\cite{Agr2} and the paper of Falbel and Gorodski~\cite{FalGor} for background on these structures.

Let~$M$ be a~smooth manifold.
A~\textit{sub-Lorentzian structure} on~$M$ is a~pair $(H,g)$, where~$H$ is a~bracket generating distribution of constant
rank on~$M$, and~$g$ is a~Lorentzian metric on~$H$.
A~triple $(M,H,g)$, where $(H,g)$ is a~sub-Lorentzian structure on~$M$, will be called a~\textit{sub-Lorentzian manifold}.

For any $q\in M$, a~vector $v\in H_{q}$ will be called \textit{horizontal}.
A~vector f\/ield~$X$ on~$M$ is horizontal if it takes values in~$H$.
We will denote the set of all local horizontal vector f\/ields by $\Gamma(H)$.
To be more precise, $X \in \Gamma (H)$ if and only if~$X$ is a~horizontal vector f\/ield def\/ined on some open subset
$U\subset M$.

A nonzero vector $v\in H_{q}$ is said to be \textit{timelike} (resp.\
\textit{spacelike}, \textit{null, nonspacelike}) if $g(v,v)<0$ (resp.\
$g(v,v)>0$, $g(v,v)=0$, $g(v,v)\leq 0$), moreover the zero vector is def\/ined to be spacelike.
Similarly, vector f\/ields are categorized analogously according to when their values lie in exactly one of the four
categories mentioned above.
An absolutely continuous curve $\gamma:[a,b]\longrightarrow M$ is called horizontal if $\dot{\gamma}(t)\in H_{\gamma(t)}$ a.e.\
on $[a,b]$.
A~horizontal curve $\gamma:[a,b]\longrightarrow M$ is timelike (spacelike, null, nonspacelike) if $\dot{\gamma}(t)$ is
timelike (spacelike, null, nonspacelike) a.e.~on $[a,b]$.

If $(H,g)$ is a~sub-Lorentzian metric on~$M$ then, as shown in~\cite{groGlob},~$H$ can be represented as a~direct sum
$H=H^{-}\oplus H^{+}$ of subdistributions such that $\rank H^{-}=1$ and the restriction of~$g$ to~$H^{-}$
(resp.\
to $H^{+}$) is negative (resp.\
positive) def\/inite.
This type of decomposition will be called \textit{a causal decomposition of}~$H$.
Now, by \textit{a time $($resp.\
space$)$ orientation} of $(M,H,g)$ we mean an orientation of the vector bundle $H^{-}\longrightarrow M$ (resp.\
$H^{+}\longrightarrow M$).
This def\/inition requires some explanations since causal decompositions are not unique.
So suppose that we are given two causal decompositions $H=H_{1}^{-}\oplus H_{1}^{+}=H_{2}^{-}\oplus H_{2}^{+}$ and that
the bundles $H_{i}^{-}\longrightarrow M$ (resp.\
$H_{i}^{+}\longrightarrow M)$, $i=1,2$, are oriented.
We say that the two orientations of $H_{1}^{-}\longrightarrow M$ and $H_{2}^{-}\longrightarrow M$ (resp.\
of $H_{1}^{+}\longrightarrow M$ and $H_{2}^{+}\longrightarrow M$) def\/ine the same time (resp.\
space) orientation of $(M,H,g)$ if around any point of~$M$ there exist local sections $X_{1}^{(i)}$ of $H_{i}^{-}\longrightarrow M$
which agree with the given orientations of $H_{i}^{-}\longrightarrow M$ (resp.\
local sections $X_{2}^{(i)},\dots,X_{k}^{(i)}$ of $H_{i}^{+}\longrightarrow M$ which agree with the given orientations of
$H_{i}^{+}\longrightarrow M$), $i=1,2$, such that $g\big(X_{1}^{(1)},X_{1}^{(2)}\big)<0$ (resp.\
$\det \big(g\big(X_{i}^{(1)},X_{j}^{(2)}\big)_{i,j=2,\dots,k}>0\big)$; by~$k$ we denote the rank of~$H$.

Since a~line bundle is orientable if and only if it is trivial, time orientability of $(M,H,g)$ is equivalent to the
existence of a~continuous timelike vector f\/ield on~$M$.
A~choice of such a~timelike f\/ield is called a~time orientation of $(M,H,g)$.

\looseness=-1
Suppose that $(M,H,g)$ is time oriented by a~vector f\/ield~$X$.
A~nonspacelike $v\in H_{q}$ will be called \textit{future} (resp.\
\textit{past}) \textit{directed} if $g(v,X(q))<0$ (resp.\
$g(v,X(q))>0$).
A~horizontal curve $\gamma:[a,b]\longrightarrow M$ is called timelike future (past) directed if $\dot{\gamma}(t)$ is
timelike future (past) directed a.e.
Similar classif\/ications can be made for other types of curves, e.g.\
nonspacelike future directed etc.
If $q_{0}\in M$ is a~point and~$U$ is a~neighborhood of $q_{0}$, then by the future timelike (nonspacelike, null)
reachable set from $q_{0}$ relative to~$U$ we mean the set of endpoints of all timelike (nonspacelike, null) future
directed curves that start from $q_{0}$ and are contained in~$U$.

Now we def\/ine a~very important notion that will play a~crucial role in the sequel.
As it is known~\cite{g6}, any sub-Lorentzian structure $(H,g)$ determines the so-called \textit{geodesic Hamiltonian}
which is def\/ined as follows.
The existence of the structure $(H,g)$ is equivalent to the existence of the f\/iber bundle morphism
$G:T^{\ast}M\longrightarrow H$ covering the identity def\/ined by $G(\lambda)=(\lambda_p |_{H_p})^\sharp$,
where $\lambda \in T_p^*M$ and $\sharp$ denotes the musical isomorphism.
In particular, if $v$, $w$ are any horizontal vectors, then $g(v,w)=\left\langle \xi,G\eta \right\rangle=\left\langle
\eta,G\xi \right\rangle$ whenever $\xi \in G^{-1}(v)$, $\eta \in G^{-1}(w)$.
The geodesic Hamiltonian is the map $h:T^{\ast}M\longrightarrow \mathbb{R}$ def\/ined~by
\begin{gather*}
h(\lambda)=\frac{1}{2}\left\langle \lambda,G\lambda \right\rangle.
\end{gather*}
If $X_{1},\dots,X_{k}$ is an orthonormal basis for $(H,g)$ with a~time orientation $X_{1}$, then
\begin{gather*}
h|_{T_{q}^{\ast}M}(\lambda)=-\frac{1}{2}\left\langle
\lambda,X_{1}(q)\right\rangle^{2}+\frac{1}{2}\sum\limits_{i=2}^{k}\left\langle \lambda,X_{i}(q)\right\rangle^{2}.
\end{gather*}
A~horizontal curve $\gamma:[a,b]\longrightarrow M$ is said to be a~\textit{Hamiltonian geodesic} if there exists
$\Gamma:[a,b]\longrightarrow T^{\ast}M$ such that $\dot{\Gamma}=\vec{h}(\Gamma)$ and $\pi (\Gamma (t))=\gamma(t)$
on $[a,b]$; by $\pi:T^{\ast}M\longrightarrow M$ we denote the canonical projections, and $\vec{h}$ is
the Hamiltonian vector f\/ield corresponding to~$h$.

Let $\gamma:[a,b]\longrightarrow M$ be a~nonspacelike curve.
The non-negative number
\begin{gather*}
L(\gamma)=\int_{a}^{b}\left\vert g(\dot{\gamma}(t), \dot{\gamma}(t))\right\vert^{1/2}dt
\end{gather*} is
called the \textit{sub-Lorentzian length of a~curve} $\gamma$.
If $U\subset M$ is an open subset, then the \textit{$($local$)$ sub-Lorentzian distance relative to} $U$
is the function $d[U]:U\times U\longrightarrow \lbrack 0,+\infty]$ def\/ined as follows: For $q_{1},q_{2}\in U$, let
$\Omega_{q_{1},q_{2}}^{nspc}(U)$ denote the set of all nonspacelike future directed curves contained in~$U$ which join
$q_{1}$ to $q_{2}$, then
\begin{gather*}
d[U](q_{1},q_{2})=\begin{cases}
\sup \left\{L(\gamma):\gamma \in \Omega_{q_{1},q_{2}}^{nspc}(U)\right\}: & \Omega_{q_{1},q_{2}}^{nspc}(U)\neq \varnothing,
\\
0: & \Omega_{q_{1},q_{2}}^{nspc}(U)=\varnothing
\end{cases}
\end{gather*}
(if $\Omega_{q,q}^{nspc}(U)$ is non-empty for a~$q\in U$ then $d[U](q,q)=\infty $).
A~nonspacelike future directed curve $\gamma:[a,b]\longrightarrow U$ is called a~$U$-\textit{maximizer} if $d[U](\gamma
(a),\gamma (b))=L(\gamma)$.
It can be proved (see~\cite{g6}) that every suf\/f\/iciently small subarc of every nonspacelike future directed Hamiltonian
geodesic is a~$U$-maximizers for suitably chosen~$U$.

Suppose now that we are given two sub-Lorentzian manifolds $(M_{i},H_{i},g_{i})$, $i=1,2$.
A~dif\/feo\-morphism $\varphi:M_{1}\longrightarrow M_{2}$ is said to be a~\textit{sub-Lorentzian isometry}, if $d\varphi
(H_{1})\subset H_{2}$ and for each $q\in M_{1}$, the mapping $d\varphi_{q}:(H_{1})_{q}\longrightarrow (H_{2})_{\varphi
(q)}$ is a~linear isometry, i.e., for every $v_{1},v_{2}\in (H_{1})_{q}$ it follows that
\begin{gather*}
g_{1}(v_{1},v_{2})=g_{2}(d\varphi_{q}(v_{1}),d\varphi_{q}(v_{2})).
\end{gather*}
Of course, any isometry maps timelike curves from $M_{1}$ to timelike curves on $M_{2}$.
The same for spacelike and null curves.
Moreover isometries preserve the sub-Lorentzian length of nonspacelike curves.
If $(M_{i},H_{i},g_{i})$, $i=1,2$, are both time- and space-oriented, then we can distinguish among all isometries
those that preserve one of the orientations or both of them.
More precisely, suppose that $\varphi:M_{1}\longrightarrow M_{2}$ is an isometry.
Let $H_{1}=H_{1}^{-}\oplus H_{1}^{+}$ be a~causal decomposition with given orientation on $H_{1}^{\pm}$.
Let $(H_{2}^{-})_{\varphi (q)}=d\varphi_{q}(H_{1}^{-})_{q}$, and $(H_{2}^{+})_{\varphi
(q)}=d\varphi_{q}(H_{1}^{+})_{q}$, $q\in M$.
Then $H_{2}=H_{2}^{-}\oplus H_{2}^{+}$ is again a~causal decomposition where the summands $H_{2}^{\pm}$ inherit the
orientation carried from~$H_{1}^{\pm}$ by~$\varphi$.
Now we say that $\varphi$ preserves time (resp.\
space) orientation if the orientations of~$H_{2}^{-}$ (resp.\
$H_{2}^{+}$) induced by $\varphi$ agrees with the time (resp.\
space) orientation of $(M_{2},H_{2},g_{2})$.
An isometry that preserves time and space orientation will be called \textit{a} $ts$\textit{-isometry}.
It is clear that any $ts$-isometry preserves Hamiltonian geodesics, maximizers, and local sub-Lorentzian distance
functions.
Notice furthermore that the set of all isometries $(M,H,g)\longrightarrow (M,H,g)$ is a~Lie group and the set of all
$ts$-isometries forms a~connected component containing the identity.

A sub-Lorentzian manifold $(M,H,g)$ is called a~\textit{contact sub-Lorentzian manifold}, if~$H$ is a~contact
distribution on~$M$.
Among sub-Lorentzian manifolds, those which are contact seem to be the easiest to study and hence well known.
Contact sub-Lorentzian manifolds are studied for instance in papers~\cite{g2, bcc2004, g7, Gro02, Grong, Huang, markin, markin2}.
The investigations go in two directions.
The f\/irst addresses global aspects, e.g., in~\cite{bcc2004, g7} the Heisenberg sub-Lorentzian metric is treated.
More precisely, the future timelike, nonspacelike and null reachable sets from a~point are computed, and a~certain
estimate on the distance function is given.
Moreover, it is shown that the future timelike conjugate locus of the origin is zero, while the future null conjugate
locus equals the union of the two null future directed Hamiltonian geodesics starting from the origin.
In turn, in~\cite{markin} and~\cite{g7} it is proved that the set reachable from the origin by future directed timelike
Hamiltonian geodesics coincides with the future timelike reachable set from the origin.
In~\cite{markin} the authors also study the set reachable by spacelike Hamiltonian geodesics and prove the uniqueness of
geodesics in the Heisenberg case.
Next, in the papers~\cite{markin, markin2} the so-called $\mathbb{H}$-type groups (i.e.~higher-dimensional analogues of
the 3D Heisenberg group) with suitable sub-Lorentzian metrics are studied, and the main emphasis is put on the problem
of connectivity by geodesics, i.e.~given two points $q_{1}$, $q_{2}$,
f\/igure out how many geodesics joining~$q_{1}$ to~$q_{2}$ exist.
A~similar problem is also dealt with in~\cite{Huang}.
On the other hand, in~\cite{Grong} the group ${\rm SL}(2,\mathbb{R})$ with the sub-Lorentzian metric is studied.
As it will become clear below, the cases of the Heisenberg group and that of ${\rm SL}(2,\mathbb{R})$ are especially
interesting for us because these are exactly the cases that arise when the invariant $\tilde{h}$ (def\/ined below) vanishes.
The other direction of studies concerns the local situation and is based on the construction (see, e.g.,~\cite{g2, Gro01, Gro02}) of local normal forms.
Such local normal forms depend on two functional parameters and permit to view general structures as perturbations of
f\/lat ones.
This fact allows to generalize some global results that hold in the f\/lat case to local results for general structures
describable by the mentioned normal forms.

As one can see, problems connected with isometric and conformal symmetry have not been examined in an explicit sense
although in broader contexts such as parabolic geometry and Cartan's equivalence, there are applicable results.
The aim of this paper is to embark on f\/illing this gap.
More precisely, f\/irst we construct invariants for contact sub-Lorentzian manifolds $(M,H,g)$ with $\dim M=3$, more or
less in the way as it is done in the contact sub-Riemannian case~-- cf.~\cite{Agr2}.
Our invariants are: a~$(1,1)$-tensor $\tilde{h}$ on~$H$ and a~smooth function $\kappa$ on~$M$.
Then, we consider in some detail the case that~$M$ is a~$3$-dimensional Lie group such that $\tilde{h}=0$.
It turns out that in such a~case~$M$ is locally either the Heisenberg group or the universal cover of
${\rm SL}(2,\mathbb{R})$.
In these two cases we describe inf\/initesimal isometries and more generally inf\/initesimal conformal transformations.

\subsection{The content of the paper}

In Section~\ref{Section2} we construct invariants for $ts$-oriented contact sub-Lorentzian metrics on 3D~manifolds.
The construction follows the ideas of~\cite{Agr2}, however the full analogy does not exist due to the special character
of indef\/inite case.
Our main invariants for a~manifold $(M,H,g)$ are: a~smooth $(1,1)$-tensor $\tilde{h}$ on~$H$ and a~smooth function $\kappa$ on~$M$.
These invariants provide necessary conditions for two contact sub-Lorentzian manifolds to be locally $ts$-isometric.
We also consider another invariant $\chi$ arising from the eigenvalues of $\tilde{h}$ which to a~lesser extent also
distinguishes the structure.
The question as to whether $\{\tilde h, \kappa \}$ is a~complete set of invariants requires deeper ana\-ly\-sis using
Cartan's theory and is deferred to a~forthcoming paper with Alexandr Medvedev~\cite{GMW}.

In Section~\ref{SubLorentzInf} we def\/ine and prove basic properties of inf\/initesimal sub-Lorentzian isometries and conformal transformations.
Then we notice that the invariant $\tilde{h}$ can be expressed in terms of the restricted Lie derivative of the
metric~$g$ in the direction of the Reeb vector f\/ield.
The immediate consequence of this latter fact is that the Reeb vector f\/ield $X_{0}$ is an inf\/initesimal isometry if and
only if $\tilde{h}$ vanishes identically.

Section~\ref{Section4} covers some other implications of certain combinations of the invariants vanishing.
In particular we demonstrate (see Proposition~\ref{smpcon}) that without any assumptions on orientation, the condition
$\chi=0$ and $\tilde{h}\neq 0$ implies the existence of line sub-bundle $L\rightarrow M$ of~$H$ on which the
metric~$g$ is equal to zero.
We then begin to focus on the condition $\tilde{h}=0$, where $\kappa$ comes to the fore.
For example, when~$M$ is a~simply connected Lie group, we show that $\tilde{h}=0$ and $\kappa=0$ implies that~$M$ is
the Heisenberg group~-- cf.~Corollary~\ref{CorHeis}, and $\tilde{h}=0$ and $\kappa \neq 0$ implies that $M$ is the
universal cover of ${\rm SL}(2,\mathbb{R})$~-- see Corollary~\ref{CorSL2}.
This contrasts with the sub-Riemannian case where a~third group, namely ${\rm SU}(2)$, also appears.

Section~\ref{Section6} is devoted to computing inf\/initesimal isometries and inf\/initesimal conformal transformations using Cartan's
equivalence method and Appendix
presents an example of an isometrically rigid sub-Lorentzian structure.

Finally, the appendix presents possible applications of our invariants to a~non-contact case.

\section{Constructing the invariants}\label{Section2}

\subsection{Preliminaries}

Let $(M,H,g)$ be a~contact sub-Lorentzian manifold, $\dim M=3$, which is supposed to be both time and space oriented or
$ts$\textit{-oriented} for short.
Since~$H$ is of rank~$2$, any causal decomposition $H=H^{-}\oplus H^{+}$ splits~$H$ into a~direct sum of line bundles.
So in this case a~space orientation is just a~continuous spacelike vector f\/ield, and consequently~$H$ admits a~global basis.
Let us f\/ix an orthonormal basis $X_{1},X_{2}$ for $(H,g)$, i.e.
\begin{gather*}
g(X_{1},X_{1})=-1,
\qquad
g(X_{1},X_{2})=0,
\qquad
g(X_{2},X_{2})=1,
\end{gather*}
where $X_{1}$ (resp.\
$X_{2}$) is a~time (resp.\
space) orientation.
From now on we will work with $ts$\textit{-invariants}, i.e.,~with invariants relative to $ts$-isometries.
However, when reading the text the reader will see that the space orientation is only an auxiliary notion here and most
of the results do not depend on it (some of them do not depend on an orientation at all).

Let $\omega$ be a~contact $1$-form such that $H=\ker \omega$.
Without loss of generality we may assume that $\omega$ is normalized so that
\begin{gather*}
d\omega (X_{1},X_{2})=\omega ([X_{2},X_{1}])=1.
\end{gather*}
Next, denote by $X_{0}$ the so-called \textit{Reeb vector field} on~$M$ which is def\/ined by
\begin{gather}
\omega (X_{0})=1, \qquad d\omega (X_{0},\cdot)=0.
\label{Reeb}
\end{gather}
It is seen that $X_{0}$ is uniquely determined by $ts$-oriented sub-Lorentzian structure.
Using~\eqref{Reeb} it is seen that the action of $\operatorname{ad}_{X_{0}}$ preserves the horizontality of vector f\/ields, i.e.
\begin{gather}
\operatorname{ad}_{X_{0}}(\Gamma (H))\subset \Gamma (H).
\label{act}
\end{gather}

Now (similarly as in~\cite{Agr2}) we introduce the structure functions.
Thanks to~\eqref{act} and~\eqref{Reeb} we have
\begin{gather}
[X_{1},X_{0}]=c_{01}^{1}X_{1}+c_{01}^{2}X_{2},
\qquad
[X_{2},X_{0}]=c_{02}^{1}X_{1}+c_{02}^{2}X_{2},
\nonumber
\\
[X_{2},X_{1}]=c_{12}^{1}X_{1}+c_{12}^{2}X_{2}+X_{0}.
\label{StrConst}
\end{gather}

Let $\nu_{0}$, $\nu_{1}$, $\nu_{2}$ be the dual basis of $1$-forms: $ \langle \nu_{i},X_{j} \rangle=\delta_{ij}$,
$i,j=0,1,2$.
Rewriting~\eqref{StrConst} in terms of $\nu_{i}$'s we have
\begin{gather}
d\nu_{0}=\nu_{1}\wedge \nu_{2},
\qquad
d\nu_{1}=c_{01}^{1}\nu_{0}\wedge \nu_{1}+c_{02}^{1}\nu_{0}\wedge \nu_{2}+c_{12}^{1}\nu_{1}\wedge \nu_{2},
\nonumber
\\
d\nu_{2}=c_{01}^{2}\nu_{0}\wedge \nu_{1}+c_{02}^{2}\nu_{0}\wedge \nu_{2}+c_{12}^{2}\nu_{1}\wedge \nu_{2}.
\label{StrConstForm}
\end{gather}
Dif\/ferentiating the f\/irst equation in~\eqref{StrConstForm} we obtain $0=d\nu_{1}\wedge \nu_{2}-\nu_{1}\wedge
d\nu_{2}=(c_{01}^{1}+c_{02}^{2})\nu_{0}\wedge \nu_{1}\wedge \nu_{2}$ from which it follows that
\begin{gather}
c_{01}^{1}+c_{02}^{2}=0.
\label{eq1}
\end{gather}

\subsection{Induced bilinear form and linear operator}

In the introduction we def\/ined the geodesic Hamiltonian~$h$ which can be written as
\begin{gather*}
h=-\frac{1}{2}h_{1}^{2}+\frac{1}{2}h_{2}^{2},
\end{gather*}
where $h_{i}(\lambda)=\langle \lambda,X_{i}\rangle$, $i=1,2$.
We also consider the function $h_{0}(\lambda)=\langle \lambda,X_{0}\rangle$ and observe that by def\/inition
both~$h$ and $h_{0}$ are invariant with respect to $ts$-oriented structure $(H,g)$.
Therefore, it is the same with their Poisson bracket $\{h,h_{0}\}$ which, when evaluated at $q\in M$, gives
a~symmetric bilinear form
\begin{gather*}
\{h,h_{0}\}_{q}: \ T_{q}^{\ast}M\times T_{q}^{\ast}M\longrightarrow \mathbb{R}.
\end{gather*}
If $\lambda \in T_{q}^{\ast}M$ then $\lambda=\sum\limits_{i=0}^{3}h_{i}(\lambda)\nu_{i}(q)$ and

\begin{Lemma}
\label{ala1}
$\{h,h_{0}\}_{q}=-c_{01}^{1}h_{1}^{2}+(c_{02}^{1}-c_{01}^{2})h_{1}h_{2}+c_{02}^{2}h_{2}^{2}$.
\end{Lemma}

\begin{proof}
The formula follows from $\{h,h_{0}\}=-h_{1}\{h_{1},h_{0}\}+h_{2}\{h_{2},h_{0}\}$,
where we substitute $\{h_{i},h_{0}\} (\lambda)=\langle \lambda, [X_{i},X_{0}] \rangle$, and then use~\eqref{StrConst}.
\end{proof}

In the assertion of Lemma~\ref{ala1} and in many other places below we write $c_{jk}^{i}$ for $c_{jk}^{i}(q)$.

It follows that $ \{h,h_{0} \}_{q}(\lambda,\cdot)=0$ whenever $\lambda \in H_{q}^{\perp}$ (by def\/inition
$H_{q}^{\perp}$ is the set of such covectors $\lambda \in T_{q}^{\ast}M$ that $ \langle \lambda,v \rangle=0$
for every $v\in H_{q}$), so in fact
\begin{gather*}
 \{h,h_{0} \}_{q}: \ T_{q}^{\ast}M/H_{q}^{\perp}\times T_{q}^{\ast}M/H_{q}^{\perp}\longrightarrow
\mathbb{R}.
\end{gather*}

Let us recall (mutually inverse) musical isomorphisms determined by the metric~$g$; these are
$^{\sharp}:H^{\ast}\longrightarrow H$ and $^{\flat}:H\longrightarrow H^{\ast}$, where by def\/inition
$(\nu_{1})^{\sharp}=-X_{1}$, $(\nu_{2})^{\sharp}=X_{2}$, $(X_{1})^{\flat}=-\nu_{1}$,
$(X_{2})^{\flat}=\nu_{2}$.
Now it is easy to see that the bundle morphism $G:T^{\ast}M\longrightarrow H$ from the introduction induces for
each~$q$ a~natural identif\/ication
\begin{gather*}
F_{q}: \ T_{q}^{\ast}M/H_{q}^{\perp}\longrightarrow H_{q},
\qquad
F([\alpha])= (\alpha_{|H_{q}} )^{\sharp},
\end{gather*}
where $[\alpha]$ stands for the class of $\alpha \in T_{q}^{\ast}M$ modulo $H_{q}^{\perp}$; more precisely,
$H_{q}^{\perp}$ is spanned by~$\nu_{0}$, and $F([\nu_{1}])=-X_{1}$, $F([\nu_{2}])=X_{2}$.
This permits us to def\/ine a~bilinear symmetric form $\bar{h}_{q}:H_{q}\times H_{q}\longrightarrow \mathbb{R}$~by
\begin{gather*}
\bar{h}_{q}(v,w)= \{h,h_{0} \}_{q}\big(F_{q}^{-1}(v),F_{q}^{-1}(w)\big) .
\end{gather*}
Its matrix in the basis $X_{1}(q)$, $X_{2}(q)$ is
\begin{gather}
\left(
\begin{matrix}
-c_{01}^{1}
&-\frac{1}{2}\big(c_{02}^{1}-c_{01}^{2}\big)
\\
-\frac{1}{2}\big(c_{02}^{1}-c_{01}^{2}\big) & c_{02}^{2}
\end{matrix}
\right) .
\label{mac}
\end{gather}
Finally we def\/ine a~linear mapping $\tilde{h}_{q}:H_{q}\longrightarrow H_{q}$ by the following formula:
\begin{gather*}
\tilde{h}_{q}(v)=\left(\bar{h}_{q}(v,\cdot)\right)^{\sharp}.
\end{gather*}
Using~\eqref{mac}, it is seen that the matrix of the operator $\tilde{h}_{q}$ in the basis $\{X_{1}(q),X_{2}(q)\}$ is
equal to
\begin{gather}
\left(
\begin{matrix}
c_{01}^{1} & \frac{1}{2}(c_{02}^{1}-c_{01}^{2})
\\
-\frac{1}{2}(c_{02}^{1}-c_{01}^{2}) & c_{02}^{2}
\end{matrix}
\right).
\label{mac2}
\end{gather}

\subsection[The $ts$-invariants]{The $\boldsymbol{ts}$-invariants}

By our construction, the eigenvalues and determinant of $\tilde{h}_{q}$ as well as $\tilde{h}_{q}$ itself, are all
invariants for the $ts$-oriented structure $(H,g)$.
Clearly $\det\tilde{h}_{q}=c_{01}^{1}c_{02}^{2}+\frac{1}{4}(c_{02}^{1}-c_{01}^{2})^{2}=-(c_{01}^{1})^{2}+\frac{1}{4}(c_{02}^{1}-c_{01}^{2})^{2}$.
Since, in view of \eqref{eq1}, the trace of $\tilde{h}_{q}$ is equal to~$0$, the eigenvalues of $\tilde{h}_{q}$ are
equal to $\pm \sqrt{-(c_{01}^{1})^{2}+\frac{1}{4}(c_{02}^{1}-c_{01}^{2})^{2}}$.
We can choose
\begin{gather*}
\chi=-\big(c_{01}^{1}\big)^{2}+\frac{1}{4}\big(c_{02}^{1}-c_{01}^{2}\big)^{2}
\end{gather*}
as a~functional $ts$-invariant for our structure.
In analogy with the sub-Riemannian case~\cite{Agr2,Agr1,Agr0}, we consider the functional
$ts$-invariant def\/ined as follows:
\begin{gather*}
\kappa=X_{2}\big(c_{12}^{1}\big)+X_{1}\big(c_{12}^{2}\big)-\big(c_{12}^{1}\big)^{2}+\big(c_{12}^{2}\big)^{2}
-\frac{1}{2}\big(c_{01}^{2}+c_{02}^{1}\big).
\end{gather*}
Unlike the sub-Riemannian case where $\chi$ and $\kappa$ play the crucial role, it is $\tilde{h}$ and $\kappa$ that
play the crucial role in the sub-Lorentzian setting.

\begin{Proposition}\label{oro3}
$\kappa$ is indeed a~$ts$-invariant.
\end{Proposition}

\subsubsection{Proof of Proposition~\ref{oro3}}

Let $X_{1}$, $X_{2}$ is an orthonormal basis with a~time orientation $X_{1}$ and a~space orientation~$X_{2}$, and let
$c_{jk}^{i}$ be structures functions determined by this basis.
Next, let $\theta=\theta (q)$ be a~smooth functions and consider an orthonormal basis $Y_{1}$, $Y_{2}$ given~by
\begin{gather}
Y_{1}=X_{1}\cosh \theta+X_{2}\sinh \theta,
\qquad
Y_{2}=X_{1}\sinh \theta+X_{2}\cosh \theta.
\label{transf1}
\end{gather}
Then $Y_{1}$ ($Y_{2}$) is a~time (space) orientation, and of course
\begin{gather*}
X_{1}=Y_{1}\cosh \theta-Y_{2}\sinh \theta,
\qquad
X_{2}=-Y_{1}\sinh \theta+Y_{2}\cosh \theta.
\end{gather*}
Let $d_{jk}^{i}$ be the structure functions determined by the basis $Y_{1}$, $Y_{2}$, i.e.
\begin{gather*}
[Y_{1},Y_{0}]=d_{01}^{1}Y_{1}+d_{01}^{2}Y_{2},
\qquad
[Y_{2},Y_{0}]=d_{02}^{1}Y_{1}+d_{02}^{2}Y_{2},
\\
[Y_{2},Y_{1}]=d_{12}^{1}Y_{1}+d_{12}^{2}Y_{2}+X_{0}.
\end{gather*}
In order to prove Proposition~\ref{oro3} we need the following lemma.

\begin{Lemma}
\label{oro4}
The following formulas hold true:
\begin{gather}
d_{02}^{1}=-X_{0}(\theta)+c_{02}^{1}\cosh^{2}\theta-c_{01}^{2}\sinh^{2}\theta+\big(c_{01}^{1}-c_{02}^{2}\big) \sinh \theta \cosh \theta,
\nonumber
\\
d_{02}^{2}=\big(c_{01}^{2}-c_{02}^{1}\big) \sinh \theta \cosh \theta+c_{02}^{2}\cosh^{2}\theta-c_{01}^{1}\sinh^{2}\theta,
\nonumber
\\
d_{01}^{1}=c_{01}^{1}\cosh^{2}\theta-c_{02}^{2}\sinh^{2}\theta+\big(c_{02}^{1}-c_{01}^{2}\big) \sinh \theta \cosh \theta,
\nonumber
\\
d_{01}^{2}=-X_{0}(\theta)+c_{01}^{2}\cosh^{2}\theta-c_{02}^{1}\sinh^{2}\theta+\big(c_{02}^{2}-c_{01}^{1}\big) \sinh \theta \cosh \theta,
\nonumber
\\
d_{12}^{1}=\big(c_{12}^{1}-X_{1}(\theta)\big) \cosh \theta-\big(X_{2}(\theta)+c_{12}^{2}\big) \sinh \theta,
\nonumber
\\
d_{12}^{2}=\big(X_{1}(\theta)-c_{12}^{1}\big) \sinh \theta+\big(X_{2}(\theta)+c_{12}^{2}\big) \cosh \theta.
\label{StrConst2}
\end{gather}
\end{Lemma}

\begin{proof}
All formulas are proved by direct calculations.
For instance, using \eqref{transf1} we write
\begin{gather*}
[Y_{2},Y_{1}]=[X_{1}\sinh \theta+X_{2}\cosh \theta,X_{1}\cosh \theta+X_{2}\sinh \theta ]
\\
\phantom{[Y_{2},Y_{1}]}
=-X_{1}(\theta)X_{1}+X_{2}(\theta)X_{2}+[X_{2},X_{1}] .
\end{gather*}
Then using~\eqref{transf1} and~\eqref{StrConst} we arrive at
\begin{gather*}
-X_{1}(\theta)X_{1}+X_{2}(\theta)X_{2}+c_{12}^{1}X_{1}+c_{12}^{2}X_{2}+X_{0}
\\
\qquad
=-X_{1}(\theta)(Y_{1}\cosh \theta-Y_{2}\sinh \theta )+\allowbreak X_{2}(\theta)(-Y_{1}\sinh \theta+Y_{2}\cosh \theta )
\\
\qquad
\phantom{=}{}
+c_{12}^{1}(Y_{1}\cosh \theta-Y_{2}\sinh \theta )+\allowbreak c_{12}^{2}(-Y_{1}\sinh \theta
+Y_{2}\cosh \theta )+X_{0},
\end{gather*}
from which the f\/ifth and sixth equations in~\eqref{StrConst2} follow.
\end{proof}

Now, using Lemma~\ref{oro4}, we see that
\begin{gather}
\frac{1}{2}\left(d_{01}^{2}+d_{02}^{1}\right)=-X_{0}(\theta)+\frac{1}{2} \left(c_{01}^{2}+c_{02}^{1}\right) .
\label{shortref1}
\end{gather}
and
\begin{gather}
\left(d_{12}^{1}\right)^{2}-\left(d_{12}^{2}\right)^{2}=\left(X_{1}(\theta)-c_{12}^{1}\right)^{2}-\left(X_{2}(\theta)+c_{12}^{2}\right)^{2}.
\label{shortref2}
\end{gather}
Finally, we compute $Y_{2}(d_{12}^{1})+Y_{1}(d_{12}^{2})$.
To this end let us write
\begin{gather}
Y_{2}(d_{12}^{1})+Y_{1}(d_{12}^{2})=I+II,
\label{shortref3}
\end{gather}
where
\begin{gather*}
I=X_{2}(c_{12}^{1})+X_{1}\left(c_{12}^{2}\right)-[X_{2},X_{1}](\theta)
\\
\phantom{I}
=X_{2}(c_{12}^{1})+X_{1}\left(c_{12}^{2}\right)-c_{12}^{1}X_{1}(\theta)-c_{12}^{2}X_{2}(\theta)-X_{0}(\theta)
\end{gather*}
and
\begin{gather*}
II=-c_{12}^{1}X_{1}(\theta)-X_{2}^{2}(\theta)-c_{12}^{2}X_{2}(\theta)+X_{1}^{2}(\theta).
\end{gather*}
Combining~\eqref{shortref1},~\eqref{shortref2} and~\eqref{shortref3} completes the proof of Proposition~\ref{oro3}.

In summary, our basic $ts$-invariants are: a~smooth function $\kappa$ on~$M$ and a~$(1,1)$ tensor $\tilde{h}$ on~$H$.

\section{Sub-Lorentzian inf\/initesimal isometries\\ and conformal transformations}\label{SubLorentzInf}

In this section $(M,H,g)$ is a~f\/ixed sub-semi-Riemannian manifold, $\rank H$ and $\dim M$ are arbitrary.

\begin{Definition}
A a~dif\/feomorphism $f:M\longrightarrow M$ is called a~conformal transformation of $(M,H,g)$ if
(i) $d_{q}f(H_{q}) \subseteq H_{f(q)}$ for every $q\in M$,
(ii) there exists a~function $\rho \in C^{\infty}(M)$, $\rho >0$, such that
\begin{gather*}
g(d_{q}f(v),d_{q}f(w))=\rho (q)g(v,w)
\end{gather*}
for every $q\in M$ and every $v,w\in H_{q}$.
If $\rho=1$ then~$f$ is an isometry of $(M,H,g)$.
\end{Definition}

Along with conformal transformations and isometries we consider their inf\/initesimal variants.

\begin{Definition}
A~vector f\/ield~$Z$ on $(M,H,g)$ is called an inf\/initesimal conformal transformation (resp.\
inf\/initesimal isometry) if its f\/low $\psi^t$ consists of conformal transformations (isometries).
\end{Definition}

Let us note a~simple lemma.

\begin{Lemma}
\label{lem1}
Let~$Z$ be a~vector field on~$M$ and denote by $\psi^t$ its flow.
Then the following conditions are equivalent:
\begin{enumerate}\itemsep=0pt
\item[$(a)$] $\operatorname{ad}_{Z}:\Gamma (H)\longrightarrow \Gamma (H)$;

\item[$(b)$] $d_{q}\psi^t:H_{q}\longrightarrow H_{\psi^t(q)}$ for every $q\in M$ and every~$t$ such that $\psi^t$ is defined around~$q$.
\end{enumerate}
\end{Lemma}

\begin{proof}
Although the result is known, we give a~proof for the sake of completeness.

(a) $\Rightarrow$ (b) Following~\cite{LiuSuss}, we f\/ix a~point~$q$ and consider a~basis $X_{1},\dots,X_{k}$ of~$H$
def\/ined on a~neighborhood~$U$ of $q$.
By our assumption, there exist smooth functions $\alpha_{ij}$, $i,j=1,\dots,k$, such that
$\operatorname{ad}_{Z}X_{i}=\sum\limits_{j=1}^{k}\alpha_{ij}X_{j}$ on~$U$ and it follows that if $v_{i}(t)=(\psi_*^t X_{i})
(q)=d\psi^t X_{i}(\psi^{-t}q)$ then
\begin{gather*}
\dot{v}_{i}(t)=\left(\psi_*^t \operatorname{ad}_{Z}X_{i}\right) (q)=\sum\limits_{j=1}^{k}\left(\alpha_{ij}\circ\psi^{-t}\right)
(q)\left(\psi_*^{t}X_{j}\right) (q)=\sum\limits_{j=1}^{k}\beta_{ij}(t)v_{j}(t),
\end{gather*}
where $\beta_{ij}(t)=(\alpha_{ij}\circ\psi^{-t}) (q)$.
For any covector $\lambda \in T_{q}^{\ast}M$ which annihilates $H_{q}$, i.e.,~$ \left\langle \lambda,v\right\rangle=0$
for every $v\in H_{q}$, we obtain a~system of linear dif\/ferential equations for the functions $w_{i}(t)=\left\langle
\lambda,v_{i}(t)\right\rangle$, $i=1,\dots,k$:
\begin{gather*}
\dot{w}_{i}(t)=\sum\limits_{j=1}^{k}\beta_{ij}(t)w_{j}(t)
\end{gather*}
with initial conditions $w_{i}(0)=0$, $i=1,\dots,k$, since $v_{i}(0)=X_{i}(q)\in H_{q}$.
Therefore $v_{i}(t)=0$ and $(\psi_*^{t}X_{i}) (q)\in H_{q}$ every~$t$ for which $v_{i}(t)$ is def\/ined,
$i=1,\dots,k$.

(b) $\Rightarrow$ (a) Take a~point~$q$, then for every~$t$ such that $\left\vert t\right\vert$ is suf\/f\/iciently small,
we have $(\psi_*^{t}X) (q)\in H_{q}$ and it follows that $(\operatorname{ad}_{Z}X)
(q)=\frac{d}{dt}|_{t=0}(\psi_*^{-t}X) (q)\in H_{q}$.
\end{proof}

Suppose now that $f:M\longrightarrow M$ is a~dif\/feomorphism such that $df(H)=H$ and let~$T$ be a~tensor of type $(0,2)$ on~$H$.
We def\/ine a~pull-back $\tilde{f}^{\ast}:\Gamma (H)\times \Gamma (H)\longrightarrow C^{\infty}(M)$~by
\begin{gather*}
\big(\tilde{f}^{\ast}T\big)_{q}(X,Y)=T_{f(q)}(d_{q}f(X),d_{q}f(Y)) ,
\end{gather*}
where $X,Y\in \Gamma (H)$ (tilde indicates that we restrict to horizontal vector f\/ields).
We can now reformulate the def\/inition of conformal transformations in a~manner consistent with semi-Riemannian geometry:

\textit{$f$ is a~conformal transformation of $(M,H,g)$  if and only if there exists a~function $\rho \in C^{\infty}(M)$,
$\rho >0$, such that $\tilde{f}^{\ast}g=\rho g$ $($if $\rho=1$, $f$ is an isometry$)$.}

Suppose that~$Z$ is a~vector f\/ield on~$M$ such that $\operatorname{ad}_{Z}:\Gamma (H)\longrightarrow \Gamma (H)$ and let $\psi^t$
denote the (local) f\/low of $Z$.
Using Lemma~\ref{lem1}, again by analogy to the classical geometry, we can def\/ine a~local operator
$\tilde{L}_{Z}T:\Gamma (H)\times \Gamma (H)\longrightarrow C^{\infty}(M)$ which will be called \textit{the restricted
Lie derivative}:
\begin{gather}
\big(\tilde{L}_{Z}T\big)(q)=\frac{d}{dt}\bigg|_{t=0}\big(\big(\widetilde{\psi^t}\big)^{\ast}T\big) (q).
\label{RestDef}
\end{gather}
It turns out that

\begin{Proposition}
A~vector field~$Z$ is an infinitesimal conformal transformation of $(M,H,g)$ if and only if the following conditions hold:
\begin{enumerate}\itemsep=0pt
\item[$(i)$] $\operatorname{ad}_{Z}:\Gamma (H)\longrightarrow \Gamma (H)$, and
\item[$(ii)$] there exists a~function $\mu \in C^{\infty}(M)$ such that $\tilde{L}_{Z}g=u g$.
\end{enumerate}
\end{Proposition}

\begin{proof}
Remembering that we use only horizontal vector f\/ields, the proof is the same as in the classical geometry.
Again $\psi^t$ is the f\/low of~$Z$.

``$\Rightarrow$'' By Lemma~\ref{lem1} we know that (i) is satisf\/ied.
If $(\widetilde{\psi^t})^{\ast}g=\rho_{t}g$, where for each~$t$ the function $\rho_{t}$ is smooth and positive, then it follows that
\begin{gather*}
\tilde{L}_{Z}\big(\big(\widetilde{\psi^t}\big)^{\ast}g\big)=\frac{d}{ds}\bigg|_{s=0} \big(\widetilde{\psi^{s}}\big)^{\ast}
\big(\widetilde{\psi^t}\big)^{\ast}g= \frac{d}{dt} \big(\widetilde{\psi^t}\big)^{\ast}g=\frac{d}{dt} (\rho_{t}g )=
\dot{\rho}_{t}g.
\end{gather*}
On the other hand, we also have that
\begin{gather*}
\tilde{L}_{Z}\big(\big(\widetilde{\psi^t}\big)^{\ast}g\big)=\tilde{L}_{Z}(\rho_{t}g)=X(\rho_{t})g+\rho_{t}\big(\tilde{L}_{Z}g\big),
\end{gather*}
and so we see that $\tilde{L}_{Z}g=\mu g$, where
\begin{gather*}
\mu=\frac{\dot{\rho}_{t}-X(\rho_{t})}{\rho_{t}}
\end{gather*}
(note that $\rho_{0}=1$).

``$\Leftarrow$'' From Lemma~\ref{lem1} we know that $d\psi^t$ preserves~$H$.
From (ii) and~\eqref{RestDef} we have
\begin{gather*}
\frac{d}{dt}\big(\widetilde{\psi^t}\big)^{\ast}g=\big(\widetilde{\psi^t}\big)^{\ast}\big(\tilde{L}_{Z}g\big)
=\big(\widetilde{\psi^t}\big)^{\ast}(\mu g)=\big(\mu \circ \psi^t\big)\big(\widetilde{\psi^t}\big)^{\ast}g,
\end{gather*}
which implies that $(\widetilde{\psi^t})^{\ast}g= \rho_{t}g$, where
\begin{gather*}
\rho_{t}(q)=\exp \int_{0}^{t}\mu (\psi^{s}(q)) ds.  \tag*{\qed}
\end{gather*}
\renewcommand{\qed}{}
\end{proof}

By direct calculation we obtain
\begin{gather}
\big(\tilde{L}_{Z}g\big)(X,Y)=Z(g(X,Y))-g(\operatorname{ad}_{Z}X,Y)-g(X,\operatorname{ad}_{Z}Y)
\label{LieMetr}
\end{gather}
for every $X,Y\in \Gamma (H)$, which gives the following two corollaries:

\begin{Corollary}
\label{wn1}
$Z$ is an infinitesimal conformal transformation of $(M,H,g)$ if and only if there exists a~function $\mu \in
C^{\infty}(M)$ such that for every $X,Y\in \Gamma (H)$
\begin{gather*}
Z(g(X,Y))=g(\operatorname{ad}_{Z}X,Y)+g(X,\operatorname{ad}_{Z}Y)+\mu g(X,Y).
\end{gather*}
\end{Corollary}

\begin{Corollary}
$Z$ is an infinitesimal isometry of $(M,H,g)$ if and only if for every $X,Y\in \Gamma (H)$
\begin{gather*}
Z(g(X,Y))=g(\operatorname{ad}_{Z}X,Y)+g(X,\operatorname{ad}_{Z}Y).
\end{gather*}
\end{Corollary}

Furthermore:

\begin{Corollary}
If~$Z$ is an infinitesimal conformal transformation or isometry of $(M,H,g)$ then for every $n\geq 2$ and every $X,Y\in
\Gamma (H)$
\begin{gather}
\sum\limits_{k=0}^{n}\binom{n}{k}g\big(\operatorname{ad}_{Z}^{k}X,\operatorname{ad}_{Z}^{n-k}Y\big)=0.
\label{r1}
\end{gather}

\begin{proof}
Fix a~point \thinspace $q\in M$.
Under the above notation, for any $n\in \mathbb{N}$ and suf\/f\/iciently small $\left\vert t\right\vert$
\begin{gather}
\rho_{t}(\psi^tq)g\big(X(\psi^tq),Y(\psi^tq)\big)=g\big(d_{\psi^tq}\psi^{-t}(X),d_{\psi^tq}\psi^{-t}(Y)\big)
\nonumber
\\
\qquad
=g\left(\sum\limits_{k=0}^{n}\frac{t^{k}}{k!}(\operatorname{ad}_{Z}^{k}X)(q),\sum\limits_{m=0}^{n}\frac{t^{m}}{m!}(\operatorname{ad}_{Z}^{m}Y)(q)\right)+o(t^{n}).
\label{r2}
\end{gather}
Using Corollary~\ref{wn1} we can remove from~\eqref{r2} terms of order $0$ and $1$ with respect to~$t$.
What we obtain is
\begin{gather*}
\sum\limits_{k=2}^{n}t^{k}\sum\limits_{i+j=k}\frac{1}{i!j!}g\big(\big(\operatorname{ad}_{Z}^{i}X\big)(q),\big(\operatorname{ad}_{Z}^{j}Y\big)(q)\big)+o(t^{n})=0 \end{gather*}
for $\left\vert t\right\vert$ suf\/f\/iciently small, which gives~\eqref{r1}.
\end{proof}
\end{Corollary}

\section{Some properties of invariants}\label{Section4}

In this section we assume all sub-Lorentzian manifolds to be $ts$-oriented.
Let us start from an obvious observation.

\begin{Proposition}
Let $(M_{i},H_{i},g_{i})$, $i=1,2$, be a~contact $3$-dimensional $ts$-oriented sub-Lorent\-zian manifolds.
Denote by $\chi_{i}$, $\kappa_{i}$, $\tilde{h}_{i}$ the corresponding objects defined by $(H_{i},g_{i})$, $i=1,2$.
If $\varphi:(M_{1},H_{1},g_{1})\longrightarrow (M_{2},H_{2},g_{2})$ is a~local $ts$-isometry, then
$\chi_{1}=\varphi^{\ast}\chi_{2}$, $\kappa_{1}=\varphi^{\ast}\kappa_{2}$, and $\tilde{h}_{1}=\varphi^{\ast}\tilde{h}_{2}$.
\end{Proposition}

Fix a~contact $3$-dimensional sub-Lorentzian manifold $(M,H,g)$.
First of all let us notice how the invariant $\tilde{h}$ can be expressed in terms of the restricted Lie derivative of
the metric~$g$ in the direction of the Reeb f\/ield.
Indeed, knowing~\eqref{LieMetr} it is clear that for every $q\in M$ and every $v,w\in H_{q}$
\begin{gather}
\bar{h}_{q}(v,w)=\frac{1}{2}\big(\tilde{L}_{X_{0}}g\big)(q)(v,w).
\label{LieInv}
\end{gather}
Such an approach allows to def\/ine higher-order invariants, namely those that correspond to the bilinear forms
\begin{gather*}
\bar{h}_{q}^{(l)}(v,w)=\frac{1}{2}\big(\tilde{L}_{X_{0}}^{l}g\big)(q)(v,w),
\qquad
l=2,3,\dots .
\end{gather*}
In this way, however, we will not obtain any formulas involving the structure functions $c_{12}^{i}$.

Using~\eqref{LieInv} we obtain the following proposition and corollary thereof.

\begin{Proposition}
\label{ala2}
The Reeb vector field $X_{0}$ is an infinitesimal isometry for $(H,g)$ if and only if $\tilde{h}_{q}=0$ for every $q\in M$.
\end{Proposition}

\begin{Corollary}
If the Reeb vector field $X_{0}$ is an infinitesimal isometry for $(H,g)$ then $\chi=0$ everywhere.
\end{Corollary}

Proposition~\ref{ala2} shows one of the ways how to produce sub-Lorentzian isometries.
This is important because we know very little examples of such maps.

Next we study the ef\/fect on the invariants when we dilate the structure.
To this end suppose that we have a~sub-Lorentzian $ts$-oriented structure $(H,g)$ which is given by an orthonormal
frame $X_{1}$, $X_{2}$ with a~time (resp.\
space) orientation $X_{1}$ (resp.~$X_{2}$).
Let $s>0$ be a~constant.
Consider the sub-Lorentzian structure $(H^{\prime},g^{\prime})$ def\/ined by assuming the frame $X_{1}^{\prime}=sX_{1}$,
$X_{2}^{\prime}=sX_{2}$ to be orthonormal with the time (resp.\
space) orientation $X_{1}^{\prime}$ (resp.\
$X_{2}^{\prime}$).
The normalized one form $\omega^{\prime}$ which def\/ines $H^{\prime}$ is given by $\omega^{\prime}=\frac{1}{s^{2}}\omega
$, i.e., $d\omega^{\prime}(X_{1}^{\prime},X_{2}^{\prime})=\omega^{\prime}([X_{2}^{\prime},X_{1}^{\prime}])=1$.
It follows that the Reeb f\/ield is now $s^{2}X_{0}$.
Then it is easy to see that~\eqref{StrConst} can be rewritten as
\begin{gather*}
\lbrack X_{1}^{\prime},X_{0}^{\prime}]={c}_{01}^{\prime 1}X_{1}^{\prime}+{c}_{01}^{\prime 2}X_{2}^{\prime},
\qquad
\lbrack X_{2}^{\prime},X_{0}^{\prime}]={c}_{02}^{\prime 1}X_{1}^{\prime}+{c}_{02}^{\prime 2}X_{2}^{\prime},
\\
\lbrack X_{2}^{\prime},X_{1}^{\prime}]={c}_{12}^{\prime 1}X_{1}^{\prime}+{c}_{12}^{\prime 2}X_{2}^{\prime}+X_{0}^{\prime},
\end{gather*}
where ${c}_{jk}^{\prime i}=sc_{jk}^{i}$.
As a~corollary we obtain

\begin{Proposition}
Let $\chi$, $\kappa$, $\tilde{h}$ $($resp.\
$\chi^{\prime}$, $\kappa^{\prime}$, $\tilde{h}^{\prime})$ be the $ts$-invariants of the sub-Lorentzian structure defined
by an orthonormal basis $X_{1}$, $X_{2}$ $($resp.\
by $X_{1}^{\prime}=s X_{1}$, $X_{2}^{\prime}=s X_{2})$.
Then
\begin{gather*}
\chi^{\prime}=s^{2}\chi , \qquad \kappa^{\prime}=s^{2}\kappa , \qquad \tilde{h}^{\prime}=s \tilde{h}.
\end{gather*}
\end{Proposition}

\subsection[The case $\protect\chi=0$, $\tilde{h}\neq 0$]{The case $\boldsymbol{\protect\chi=0}$, $\boldsymbol{\tilde{h}\neq0}$} \label{section4.1}

Next let us assume that $\chi (q)=0$ but $\tilde{h}_{q}\neq 0$ (i.e.~$c_{01}^{1}\neq 0$) everywhere.
As we shall see we are given an additional structure in this case.
Indeed, the correspondence $q\longrightarrow \ker \tilde{h}_{q}$ def\/ines an invariantly given f\/ield of directions.
We can distinguish two cases: (i) $c_{01}^{1}=\frac{1}{2}(c_{02}^{1}-c_{01}^{2})$, and (ii)
$c_{01}^{1}=-\frac{1}{2}(c_{02}^{1}-c_{01}^{2})$.
In the f\/irst case the matrix of $\tilde{h}_{q}$ is of the form
\begin{gather*}
\left(
\begin{matrix}
c_{01}^{1}
& c_{01}^{1}
\\
-c_{01}^{1} &-c_{01}^{1}
\end{matrix}
\right)
\end{gather*}
and $\ker \tilde{h}_{q}$ is spanned by $X_{1}(q)-X_{2}(q)$ for each~$q$.
In the second case the matrix of $\tilde{h}_{q}$ is equal~to
\begin{gather*}
\left(
\begin{matrix}
c_{01}^{1}
&-c_{01}^{1}
\\
c_{01}^{1} &-c_{01}^{1}
\end{matrix}
\right)
\end{gather*}
and $\ker \tilde{h}_{q}$ is spanned by $X_{1}(q)+X_{2}(q)$.
Thus in the considered case there exists a~line sub-bundle $L\longrightarrow M$ of~$H$ on which~$g$ is equal to zero.
Of course this result is trivial under assumption on $ts$-orientation because then~$H$ admits a~global orthonormal basis
$X_{1}$, $X_{2}$ and we have in fact two such subbundles, namely $\Span\{X_{1}+X_{2}\}$ and $\Span\{X_{1}-X_{2}\}$.
What is interesting here is that the condition $\chi=0$, $\tilde{h}\neq 0$ does not depend on the assumption on
orientation.
Indeed, notice that if we change a~time (resp.\
space) orientation keeping space (resp.\
time) one then $\tilde{h}$ is multiplied by $-1$ (because so is $X_{0}$).
Moreover the condition $\chi=0$, $\tilde{h}\neq 0$ means that $\tilde{h}$ is a~non-zero map with vanishing eigenvalues,
the fact being independent of possible multiplication by $-1$.
Therefore the condition $\chi=0$, $\tilde{h}\neq 0$ makes sense even for an unoriented contact sub-Lorentzian
structures.
In this way we are led to the following proposition.

\begin{Proposition}
\label{smpcon}
Suppose that $(M,H,g)$ is a~contact sub-Lorentzian manifold $($we don't make any assumptions on orientation$)$.
If $\chi (q)=0$ and $\tilde{h}_{q}\neq 0$ for every $q\in M$ then there exists a~line sub-bundle $L\longrightarrow M$
of~$H$ on which~$g$ is equal to zero.
\end{Proposition}

\begin{proof}
Fix an arbitrary point $q\in M$.
Let $Y_{1}$, $Y_{2}$ be an orthonormal basis for $(H,g)$ def\/ined on a~neighborhood~$U$ of~$q$, where $Y_{1}$ is timelike and $Y_{2}$ is spacelike.
Supposing $Y_{1}$ (resp.\
$Y_{2}$) to be a~time (resp.\
space) orientation we can apply the above construction of $ts$-invariants obtaining the corresponding objects
$\chi_{U}$ and $\tilde{h}_{U}$.
By our assumption and the above remark $\chi_{U}=0$, $\tilde{h}_{U}\neq 0$ on~$U$, and we get an invariantly def\/ined
line sub-bundle $L_{U}\longrightarrow U$: $U\ni q\longrightarrow \ker (\tilde{h}_{U})_{q}=:L_{U}(q)$.
We repeat the same construction around any point $q\in M$, which results in the family $\left\{L_{U}\longrightarrow
U\right\}_{U\subset M}$ of line sub-bundles, indexed by elements~$U$ of an open covering of~$M$.
By construction $L_{U}(q)=L_{U^{\prime}}(q)$ for any $q\in U\cap U^{\prime}$.
\end{proof}

Let us note that if~$M$ is simply connected, then the assertion of Proposition~\ref{smpcon} holds true no matter the
values of $\chi$ and $\tilde{h}$ are, because in this case the metric $(H,g)$ admits a~global orthonormal frame,
see~\cite{groGlob}.

\subsection[The case $\tilde{h}=0$]{The case $\boldsymbol{\tilde{h}=0}$}\label{section4.2}

We begin with the following proposition which is clear because $X_{0}$ is an inf\/initesimal isometry.

\begin{Proposition}
If $\tilde{h}=0$ then $X_{0}(\kappa)=0$, i.e., $\kappa$ is constant along the trajectories of $X_{0}$.
\end{Proposition}

By~\eqref{mac2} the assumption $\tilde{h}=0$ implies
\begin{gather*}
c_{01}^{1}=c_{02}^{2}=0,
\qquad
c_{02}^{1}=c_{01}^{2}.
\end{gather*}
We will write $c=c_{02}^{1}=c_{01}^{2}$.
Now~\eqref{StrConst} takes the form
\begin{gather}
[X_{1},X_{0}]=cX_{2},
\qquad
[X_{2},X_{0}]=cX_{1},
\qquad
[X_{2},X_{1}]=c_{12}^{1}X_{1}+c_{12}^{2}X_{2}+X_{0}.
\label{StrConst3}
\end{gather}
Rewriting as above~\eqref{StrConst3} in terms of the dual forms $\nu_{i}$ we arrive at
\begin{gather}
d\nu_{0}=\nu_{1}\wedge \nu_{2},
\qquad
d\nu_{1}=c\nu_{0}\wedge \nu_{2}+c_{12}^{1}\nu_{1}\wedge \nu_{2},
\qquad
d\nu_{2}=c\nu_{0}\wedge \nu_{1}+c_{12}^{2}\nu_{1}\wedge \nu_{2}.
\label{StrConst3form}
\end{gather}

\begin{Lemma}
\label{oro7}
The following identities hold
\begin{gather*}
-X_{1}(c)-cc_{12}^{2}+X_{0}\big(c_{12}^{1}\big)=0,
\qquad
X_{2}(c)-cc_{12}^{1}+X_{0}\big(c_{12}^{2}\big)=0.
\end{gather*}
\end{Lemma}

\begin{proof}
The lemma is obtained upon applying the exterior dif\/ferential to both sides of the second and the third equation
in~\eqref{StrConst3form}.
\end{proof}

Our next aim, which will be achieved in the next subsection, is to f\/ind a~hyperbolic rotation of our frame $X_{1}$,
$X_{2}$ so that~\eqref{StrConst3} signif\/icantly simplif\/ies.
More precisely we want to kill the terms $c_{12}^{i}$, $i=1,2$.
To this end let us introduce the following $1$-form
\begin{gather}
\eta=(\kappa+c)\nu_{0}+c_{12}^{1}\nu_{1}-c_{12}^{2}\nu_{2}
\label{oneform}
\end{gather}
whose signif\/icance will become evident below.

\begin{Proposition}
\label{stwrdz1}
$d\eta=d\kappa \wedge \nu_{0}$.
\end{Proposition}

\begin{proof}
Computations give
\begin{gather*}
d\eta=d\kappa \wedge \nu_{0}+\big({-}X_{1}(c)-cc_{12}^{2}+X_{0}\big(c_{12}^{1}\big)\big) \nu_{0}\wedge \nu_{1}
+\big({-}X_{2}(c)+cc_{12}^{1}-X_{0}\big(c_{12}^{2}\big)\big) \nu_{0}\wedge \nu_{2}
\\
\phantom{d\eta=}{}
+\big(\kappa+c-X_{2}\big(c_{12}^{1}\big)-X_{1}\big(c_{12}^{2}\big)+\big(c_{12}^{1}\big)^{2}-\big(c_{12}^{2}\big)^{2} \big) \nu_{1}\wedge
\nu_{2}.
\end{gather*}
To end the proof we use Lemma~\ref{oro7} and the def\/inition of $\kappa$.
\end{proof}

\subsection{The simply-connected Lie group case}\label{section4.3}

Suppose that our contact sub-Lorentzian manifold $(M,H,g)$ is such that~$M$ is a~simply-connected Lie group and~$H$,~$g$
are left-invariant; this means that left translations of~$M$ are sub-Lorentzian isometries (note that any left-invariant
bracket generating distribution on a~$3$-dimensional Lie group is necessarily contact).
In such a~case, clearly, $\chi$~and~$\kappa$ are constant.
We also remark that unlike the general situation, the assumption on $ts$-orientation is no longer restrictive since
groups are parallelizable manifolds.
As above, assume that $\tilde{h}=0$ everywhere.

Recalling our aim formulated in the previous subsection we prove the following lemma.

\begin{Lemma}
There exists a~smooth function $\theta:M\longrightarrow \mathbb{R}$ such that $X_{1}(\theta)=c_{12}^{1}$,
$X_{2}(\theta)=-c_{12}^{2}$.
\end{Lemma}

\begin{proof}
Suppose that such a~function $\theta$ exists.
Then
\begin{gather*}
X_{0}(\theta)=  [X_{2},X_{1} ] (\theta)-c_{12}^{1}X_{1}(\theta)-c_{12}^{2}X_{2}(\theta)
\\
\phantom{X_{0}(\theta)}{}
= X_{2}\big(c_{12}^{1}\big)+X_{1}\big(c_{12}^{2}\big)-\big(c_{12}^{1}\big)^{2}+\big(c_{12}^{2}\big)^{2}= \kappa+c,
\end{gather*}
and it follows that
\begin{gather*}
d\theta=X_{0}(\theta)\nu_{0}+X_{1}(\theta)\nu_{1}+X_{2}(\theta)\nu_{2}=(\kappa
+c)\nu_{0}+c_{12}^{1}\nu_{1}-c_{12}^{2}\nu_{2}=\eta,
\end{gather*}
where $\eta$ is def\/ined by~\eqref{oneform}.
Thus to prove the existence of $\theta$ it is enough to show that $\eta$ is exact.
Since~$M$ is simply-connected we must show that $d\eta=0$.
This is however clear by Proposition~\ref{stwrdz1} and the fact that $\kappa$ is a~constant.
\end{proof}

Now we apply to our frame $X_{1}$, $X_{2}$, the hyperbolic rotation by the angle $\theta$ specif\/ied above.
As a~result, the frame $Y_{1}$, $Y_{2}$ given by~\eqref{transf1} satisf\/ies

\begin{Proposition}
\label{propliealg}
\begin{gather*}
[Y_{1},X_{0}]=-\kappa Y_{2},
\qquad
[Y_{2},X_{0}]=-\kappa Y_{1},
\qquad
[Y_{2},Y_{1}]= X_{0}.
\end{gather*}

\begin{proof}
It follows directly from facts proved in Subsections~\ref{section4.2},~\ref{section4.3}, from \eqref{StrConst3}, and from Lemma~\ref{oro4}.
\end{proof}
\end{Proposition}

We observe here the dif\/ference between the sub-Riemannian case where the brackets have the form
\begin{gather*}
[Y_{1},X_{0}]=\kappa Y_{2},
\qquad
[Y_{2},X_{0}]=-\kappa Y_{1},
\qquad
[Y_{2},Y_{1}]= X_{0}.
\end{gather*}
We remark that in the sub-Riemannian case Agrachev and Barilari~\cite{Agr2} obtain the following results: $\kappa=0$
implies~$M$ is isometric to the Heisenberg group, $\kappa>0$ implies~$M$ is isometric to ${\rm SU}_2$ and $\kappa<0$ implies~$M$
is isometric to the universal cover of ${\rm SL}_2$.
In the sub-Lorentzian case~${\rm SU}_2$ does not arise.
To be precise we have the following corollaries of Proposition~\ref{propliealg}.

\begin{Corollary}
\label{CorHeis}
If~$M$ is a~simply-connected Lie group such that $\tilde{h}$ and $\kappa$ vanish identically, then~$M$ is isometric to
the Heisenberg group.
\end{Corollary}

\begin{Corollary}
\label{CorSL2}
If~$M$ is a~simply-connected Lie group such that $\tilde{h}$ vanishes and $\kappa \neq 0$, then it is isometric to
a~sub-Lorentzian structure on $\widetilde{\rm SL}_{2}(\mathbb{R})$ induced by the Killing form.
\end{Corollary}

Before proving Corollary~\ref{CorSL2} let us recall some basic facts about the Killing form and Cartan decompositions.
For any Lie algebra $\mathfrak{g}$ the Killing form is the symmetric bilinear form def\/ined~by
$K(X,Y)=\mathrm{Trace}(\mathrm{ad}_X \mathrm{ad}_Y)$.
The Killing form has the following invariance properties:
\begin{enumerate}\itemsep=0pt
\item[1)] $K([X,Y],Z)=K(X,[Y,Z])$,
\item[2)] $K(T(X), T(Y))= K(X,Y)$ for all $T \in \mathrm{Aut}(\mathfrak{g})$.
\end{enumerate}

If $\mathfrak{g}$ is simple then any symmetric bilinear form satisfying the f\/irst invariance condition is a~scalar
multiple of the Killing form and Cartan's criterion states that a~Lie algebra is semisimple if and only if the Killing
form is non-degenerate.

A Cartan involution is any element $\Theta \in \mathrm{Aut}(\mathfrak{g})$ such that $\Theta^2=I$ and
\begin{gather*}
\langle X, Y \rangle_\Theta=-K(X,\Theta(Y))
\end{gather*}
is positive def\/inite.
Corresponding with~$\Theta$ we have a~Cartan decomposition $\mathfrak{g}=\mathfrak{t} \oplus \mathfrak{p}$,
where $\mathfrak{t}$ and $\mathfrak{p}$ are the eigenspaces corresponding with the eigenvalues $1$ and $-1$ respectively.
Since $\Theta$ is an automorphism, it follows that $[\mathfrak{t}, \mathfrak{t}] \subseteq \mathfrak{t}$,
$[\mathfrak{t}, \mathfrak{p}] \subseteq \mathfrak{p}$ and $[\mathfrak{p}, \mathfrak{p}] \subseteq \mathfrak{t}$.
Moreover, the Killing form is negative def\/inite on $\mathfrak{t}$ and positive def\/inite on $\mathfrak{p}$.

The standard Cartan involution on $\mathfrak{sl}_2$ is given by $\Theta(A)=-A^T$.
In this case we have that $\mathfrak{t}= \mathrm{span} \{f_1\}$ and $\mathfrak{p}= \mathrm{span} \{f_2,f_0\}$, where
\begin{gather*}
f_0=\frac{1}{2}\left(
\begin{matrix}
-1 & 0
\\
0 & 1
\end{matrix}
\right),
\qquad
f_1= \frac{1}{2} \left(
\begin{matrix}
0 & 1
\\
-1 & 0
\end{matrix}
\right),
\qquad
f_2= \frac{1}{2} \left(
\begin{matrix}
0 & 1
\\
1 & 0
\end{matrix}
\right),
\end{gather*}
and the Lie brackets are
\begin{gather*}
[f_2,f_1]=f_0,
\qquad
[f_1,f_0]=f_2,
\qquad
[f_2,f_0]=f_1.
\end{gather*}
The Killing form for $\mathfrak{sl}_2$ is given by $K(A,B)=4Tr(AB)$ and so the bilinear form $B(A, B)
=\frac{1}{2}K(A,B)$ satisf\/ies
\begin{alignat*}{4}
& B(f_0,f_0)=1,
\qquad &&
B(f_1,f_1)=-1,
\qquad&&
B(f_2,f_2)=1,&
\\
& B(f_0,f_1)=0,
\qquad&&
B(f_0,f_2)=0,
\qquad&&
B(f_1,f_2)=0.&
\end{alignat*}
Thus we have two choices: 1)~$\mathcal{H}_e= \mathrm{span} \{f_1,f_2\}$ or 2)~$\mathcal{H}_e= \mathrm{span} \{f_1,f_0\}$.
In each case, by left translation, we obtain left-invariant sub-Lorentzian structures on $\widetilde{\rm SL}_{2}(\mathbb{R})$ satisfying $\tilde h=0$.
An isometry between these two structures is induced by Lie algebra automorphism~$T$, where $Tf_0=f_2$, $Tf_1=-f_1$ and $Tf_2=f_0$.

\begin{proof}[Proof of Corollary~\ref{CorSL2}] First we observe that the matrices
\begin{gather*}
e_0=\frac{1}{2}\left(
\begin{matrix}
\kappa & 0
\\
0 &-\kappa
\end{matrix}
\right),
\qquad
e_1= \frac{1}{2}\left(
\begin{matrix}
0 & 1
\\
\kappa & 0
\end{matrix}
\right),
\qquad
e_2=\frac{1}{2}\left(
\begin{matrix}
0 & 1
\\
-\kappa & 0
\end{matrix}
\right),
\end{gather*}
form a~basis of $\mathfrak{sl}_2$ and the bracket relations are $[e_2,e_1]=e_0$, $[e_1,e_0]=-\kappa e_2$,
$[e_2,e_0]=-\kappa e_1$.
Furthermore $K(e_1,e_1)=2 \kappa$, $K(e_2,e_2)=-2 \kappa$, and $K(e_1,e_2)=0$.
Since we assume the sub-Lorentzian structure on~$M$ is left-invariant, the metric must be the left translation of the
metric $B(A, B)=-\frac{1}{2\kappa}K(A,B)$ on $T_eM=\mathfrak{sl}_2$.
Since~$M$ is simply connected it must be the universal cover of ${\rm SL}_2(\mathbb{R})$.
\end{proof}

We remark that in general
\begin{alignat*}{4}
& K(e_0,e_0)=2\kappa^2,
\qquad&&
K(e_1,e_1)=2 \kappa,
\qquad&&
K(e_2,e_2)=-2 \kappa,&
\\
& K(e_0,e_1)=0,
\qquad&&
K(e_0,e_2)=0,
\qquad&&
K(e_1,e_2)=0,&
\end{alignat*}
and so the corresponding Cartan involution is given~by
\begin{gather*}
\Theta(e_0)=-e_0,
\qquad
\Theta(e_1)=-e_1,
\qquad
\Theta(e_2)=e_2.
\end{gather*}
Hence $\mathfrak{t}=\mathrm{span}\{e_2\}$ and $\mathfrak{p}=\mathrm{span} \{e_1, e_0\}$.
If $|\kappa|\ne 1$, then since $K(e_0,e_0)=2\kappa^2$, the only choice we have for a~sub-Lorentzian structure induced~by
the Killing form is $\mathcal{H}_e= \mathrm{span} \{e_1,e_2\}$.

The null lines in $\mathfrak{sl}_2$ are $\mathrm{span}\{e_1-e_2\}$ and $\mathrm{span}\{e_1+e_2\}$.
Furthermore if we set
\begin{gather*}
n_0=e_3,
\qquad
n_1=\frac{1}{\sqrt{2}}(e_1-e_2),
\qquad
n_2=\frac{1}{\sqrt{2}} (e_1+e_2)
\end{gather*}
then
\begin{gather*}
[n_2,n_1]=n_0,
\qquad
[n_1,n_0]= \kappa n_1,
\qquad
[n_2,n_0]=-\kappa n_2.
\end{gather*}
If we set $\mathcal{H}_e= \mathrm{span} \{n_1, n_2\}$ and def\/ine
\begin{gather*}
B(n_1, n_1)=-1,
\qquad
B(n_2, n_2)=1,
\qquad
B(n_1, n_2)=0,
\end{gather*}
then the induced left-invariant structure on $\widetilde{\rm SL}_{2}(\mathbb{R})$ is isometrically distinct from the cases
above, indeed
\begin{gather*}
\tilde h= \kappa\left(
\begin{matrix}
1 & 0
\\
0 &-1
\end{matrix}
\right).
\end{gather*}

\section[Inf\/initesimal sub-Lorentzian transformations on groups with $\tilde h=0$]{Inf\/initesimal sub-Lorentzian transformations\\
on groups with $\boldsymbol{\tilde h=0}$}

\subsection{Introduction}

In this section we determine the conformal and isometry groups for the Heisenberg group and the universal cover of
${\rm SL}_2$.
In particular we will see that in both cases the local inf\/initesimal conformal transformations are given~by
$\mathfrak{sl}_3$.

In the context of this paper it would be natural to construct the vector f\/ields using the criteria developed in
Section~\ref{SubLorentzInf}, however this leads to complicated systems of PDEs
which we cannot provide explicit proof concerning solutions.
Instead we apply Cartan's equivalence method which leads to the general solution without having to solve PDEs.

In work in preparation with Alexandr Medvedev~\cite{GMW} we explore further the application of the Cartan approach to
sub-Lorentzian geometry.
In particular the invariants discussed here appear in a~much more systematic manner and a~complete description of all
left-invariant sub-Lorentzian structures on $3$-dimensional Lie groups will be given.

We also remark that on the subject of conformal classif\/ication of left-invariant sub-Rieman\-nian structures on $3$-dimensional Lie groups there is the very recent arXiv paper of Boarotto~\cite{Boar}.

\subsection{Preliminaries}

The assumption $\tilde h=0$ implies that the Lie algebra of~$M$ has the form given in Proposition~\ref{propliealg}.
So as to make the notation a~little less confusing in the calculations that ensue, we rewrite the Lie algebra in the
following form
\begin{gather*}
[X_1,X_3]= \kappa X_2,
\qquad
[X_2,X_3]= \kappa X_1,
\qquad
[X_1,X_2]= X_3.
\end{gather*}
We denote the ordered dual frame by $\theta=\{\theta_1,\theta_2,\theta_3\}$ and observe that by Cartan's formula the
structure equations of this coframe are
\begin{gather*}
d\theta_1= \kappa \theta_3 \wedge \theta_2,
\qquad
d\theta_2= \kappa \theta_3 \wedge \theta_1,
\qquad
d\theta_0= \theta_2 \wedge \theta_1.
\end{gather*}
Moreover the sub-Lorentzian metric has the form $\theta_2 \odot \theta_2-\theta_1 \odot \theta_1$.
It follows that the subgroup of $G \subset GL(3)$ which acts on~$\theta$ and leaves the metric conformally invariant
modulo terms of the form $\eta \odot \theta_3$, consists of matrices of the form
\begin{gather*}
\left(
\begin{matrix}
e^r \cosh(t) & e^r \sinh(t) & a
\\
e^r \sinh(t) & e^r \cosh(t) & b
\\
0 & 0 & e^{2r}
\end{matrix}
\right).
\end{gather*}
Hence a~local dif\/feomorphism $f: M \to M$ is conformal if and only if
\begin{gather}
f^* \theta= g_f \theta,
\label{confinveq}
\end{gather}
where $g_f:M \to G$.
Thus the conformal symmetry problem is exactly to f\/ind all local dif\/feomorphisms which satisfy~\eqref{confinveq} which
is precisely a~Cartan equivalence problem.

\subsection{Overview of Cartan's algorithm}

For the details we refer the reader to~\cite{Olver1}.
The f\/irst step in Cartan's algorithm is to pass to the equation $f^* d \theta= d (g_f \theta)$ and lift the problem to
$M \times G$.
In this context the lift of~$\theta$ is given by the partial coframe $\Theta \subset T^*(M \times G)$, where $\Theta_{(p,g)}=g\theta_p$.
The lift of~$f$ is def\/ined by $\tilde f(p, g)=(f(p),g g_f(p)^{-1})$ and it follows that
\begin{gather}
\tilde f^* \Theta= \Theta
\qquad
\text{and}
\qquad
\tilde f^* d \Theta= d\Theta.
\label{confinveq2}
\end{gather}
If $\Pi \subset T(M \times G)^*$ is any subset complementary to~$\Theta$ then the structure equations for the lifted
partial coframe take the form:
\begin{gather*}
d \Theta= \Pi\wedge \Theta+T \Theta \wedge \Theta
\end{gather*}
and~$T$ is referred to as the torsion.
It follows from~\eqref{confinveq2} that if~$f$ is a~local conformal dif\/feo\-mor\-phism then $T \circ \tilde f=T$, however
in general this equation is not the end of the story.
Roughly speaking all the higher-order coframe derivatives of~$T$ must also be invariant under composition with $\tilde f$,
see~\cite[Theorem 14.24]{Olver1} for a~precise statement.
Fortunately in our case this will not be an issue.

The main idea in the Cartan algorithm is to exploit the freedom of choice of~$\Pi$ so as to minimise the torsion~$T$.
There are three processes involved in minimising torsion, namely group reduction, absorption and prolongation.
A~reduction of the structure group can be carried out when the condition $f^* \theta= g_f \theta$ for some $g_f:M \to
G$ implies that $g_f:M \to G^{\prime}\subset G$.
Reductions reveal themselves as coef\/f\/icients in $T$ depending only on group parameters.
Such coef\/f\/icients can be set to a~convenient constant value as long as invertibilty is not violated.
An example of such a~phenomena that the reader may be familiar with occurs in the contact equivalence problem in
dimension~$3$, i.e., if initially we assume that $g=(g_{ij})$ and $g_{31}=g_{32}=0$ then it transpires that $g_{33}
=g_{11}g_{22}-g_{12}g_{21}$.

Absorption utilises the fact, that pointwise, each element of~$\Pi$ is a~Maurer-Cartan form on~$G$ plus a~linear
combination of the $\Theta_i$.
We choose the coef\/f\/icients of the $\Theta_i$ so that~$T$ has as many zero coef\/f\/icients as possible.

In summary the element $T_{ij}^k \in T$ is the coef\/f\/icient of $\Theta^i \wedge \Theta^j$ in the expression for $d
\Theta^k$.
If $T_{ij}^k$ is independent of absorption parameters then it gives a~group reduction by setting it equal to~$\pm 1$ or~$0$ and the choice is made so as not to violate invertibilty.
If $T_{ij}^k$ is dependent on absorption parameters then we solve one of these parameters so that~$T_{ij}^k=0$.
The aim of the algorithm is to reduce the group to $\{{\rm Id}\}$ through a~sequence of reduction and absorption cycles.

If after the f\/irst reduction and absorption we get $G={\rm Id}$ then the resulting~$T$ consists of the basic invariants for
the equivalence problem.
Otherwise not all group parameters are normalised and absorption parameters may remain undetermined with no torsion
coef\/f\/icients available to normalise them.
In this case the problem must be prolonged which means that the free absorption parameters are understood as the groups
parameters for a~structure group $G^{(1)}$ associated with a~new equivalence problem on $M \times G$.
Specif\/ically we write $\Pi= \varpi+F \Theta$, where~$F$ consists of the free absorption parameters, and consider the
equivalence problem for the partial coframe $\Theta \cup \varpi \subset T^*(M \times G)$ with structure group $G^{(1)}$
consisting of matrices of the form
\begin{gather*}
\left(
\begin{matrix}
I & 0
\\
F & I
\end{matrix}
\right).
\end{gather*}
In this context $\Theta \cup \Pi \subset T^*(M \times G \times G^{(1)})$ is a~lift of $\Theta \cup \varpi$ and we repeat
the procedure: augment, reduce and absorb \dots, until eventually we get an equivalence problem where the structure group
reduces to the identity and all absorption coef\/f\/icients are determined.
Of course it can happen that the process will not lead to such a~situation, but when it does, the result is the
structure equations of a~certain Cartan connection on~$M$.

\subsection{Calculations}

To begin we lift and def\/ine one forms $\Theta_i$ on $M \times G$ by setting
\begin{gather*}
\left(
\begin{matrix}
\Theta_1
\\
\Theta_2
\\
\Theta_3
\end{matrix}
\right)=\left(
\begin{matrix}
e^r \cosh(t) & e^r \sinh(t) & a
\\
e^r \sinh(t) & e^r \cosh(t) & b
\\
0 & 0 & e^{2r}
\end{matrix}
\right)\left(
\begin{matrix}
\theta_1
\\
\theta_2
\\
\theta_3
\end{matrix}
\right)
\end{gather*}
and augment the set $\{\Theta_1,\Theta_2,\Theta_3\} \subset T^*(M \times G)$ with the following forms:
\begin{gather*}
\Pi_1= \alpha_1+\frac{1}{2} b e^{-2r} \Theta_1-\frac{1}{2} a e^{-2r} \Theta_2-\big(B_1+ab e^{-4r}\big)\Theta_3,
\\
\Pi_2= \alpha_2-\frac{3}{2} a e^{-2r} \Theta_1+\frac{3}{2} b e^{-2r} \Theta_2+\big(a^2 e^{-4r}+\kappa e^{-2r}-B_2\big)\Theta_3,
\\
\Pi_3= \alpha_3- B_1 \Theta_1-B_2 \Theta_2-B_3 \Theta_3,
\\
\Pi_4=\alpha_4+\big(\big(a^2+b^2\big) e^{-4 r}-B_2\big) \Theta_1-\big(2ab e^{-4 r}+B_1\big)\Theta_2-B_4 \Theta_3.
\end{gather*}

The coef\/f\/icients of the $\Theta_i$ in $\Pi_j$ are determined by absorbing torsion and the $\alpha_j$ are the
Maurer--Cartan forms:
\begin{gather*}
\alpha_1= dr, \qquad \alpha_2= d t, \qquad \alpha_3= (da-a dr-b d t)e^{-2r}, \qquad \alpha_4= (db-b dr-a d t)e^{-2r}.
\end{gather*}
The coef\/f\/icients $B_1,\dots,B_4$ are undetermined parameters from absorption and so a~prolongation is required.
We write
\begin{gather*}
\Pi_1= \varpi_1-B_1 \Theta_3,
\qquad
\Pi_2= \varpi_2-B_2 \Theta_3,
\qquad
\Pi_3= \varpi_3-B_1 \Theta_1-B_2 \Theta_2-B_3 \Theta_3,
\\
\Pi_4= \varpi_4-B_2\Theta_1-B_1 \Theta_2-B_4 \Theta_3.
\end{gather*}
and consider the equivalence problem $M \times G$ given by the ordered basis
\begin{gather}
\{\Theta_1,\Theta_2,\Theta_3, \varpi_1, \varpi_2, \varpi_3, \varpi_4\}
\label{liftedlifted}
\end{gather}
with structure group $G^{(1)}$ consisting of matrices of the form
\begin{gather*}
\left(
\begin{matrix}
I & 0
\\
R & I
\end{matrix}
\right),
\qquad
\text{where}
\qquad
R=\left(
\begin{matrix}
0 & 0 &-B_1
\\
0 & 0 &-B_2
\\
-B_1 &-B_2 &-B_3
\\
-B_2 &-B_1 &-B_4
\end{matrix}
\right).
\end{gather*}

The ordered basis $\{\Theta_1,\Theta_2,\Theta_3, \Pi_1, \Pi_2, \Pi_3, \Pi_4\}$ is now viewed as the lift
of~\eqref{liftedlifted} to the $11$-dimensional manifold $M \times G \times G^{(1)}$ and again is augmented by forms
$\{\Omega_1,\Omega_2,\Omega_3, \Omega_4\} \subset T^*(M \times G \times G^{(1)})$.
We get the following reductions of the structure group $G^{(1)}$:
\begin{gather*}
B_2= \frac{1}{4} a^2 e^{-4r}+\frac{3}{4} b^2 e^{-4r}+\frac{1}{4} e^{-2r} \kappa,
\qquad
B_3= \frac{1}{2} b\big(\big(a^2-b^2\big) e^{-6r}+e^{-4r} \kappa\big),
\\
B_4= \frac{1}{2} a\big(\big(a^2-b^2\big) e^{-6r}+e^{-4r} \kappa\big)
\end{gather*}
and so $M \times G \times G^{(1)}$ becomes an $8$-dimensional manifold and we only need $\Omega=\Omega_4$ to augment.

Finally after absorption we arrive at the structure equations:
\begin{gather}
d \Theta_1= \Pi_1 \wedge \Theta_1+\Pi_2 \wedge \Theta_2+\Pi_3 \wedge \Theta_3,
\nonumber
\\
d \Theta_2= \Pi_1 \wedge \Theta_2+\Pi_2 \wedge \Theta_1+\Pi_4 \wedge \Theta_3,
\nonumber
\\
d \Theta_3= 2 \Pi_1 \wedge \Theta_3-\Theta_1 \wedge \Theta_2,
\nonumber
\\
d \Pi_1= \frac{1}{2} \Pi_4 \wedge \Theta_1-\frac{1}{2} \Pi_3 \wedge \Theta_2-\Omega \wedge \Theta_3,
\nonumber
\\
d \Pi_2= \frac{3}{2} \Pi_4 \wedge \Theta_2-\frac{3}{2} \Pi_3 \wedge \Theta_1,
\nonumber
\\
d \Pi_3= \Pi_3 \wedge \Pi_1-\Pi_4 \wedge \Pi_2-\Omega \wedge \Theta_1,
\nonumber
\\
d \Pi_4= \Pi_4 \wedge \Pi_1-\Pi_3 \wedge \Pi_2-\Omega \wedge \Theta_2,
\nonumber
\\
d \Omega= 2\Omega \wedge \Pi_1-\Pi_4 \wedge \Pi_3.
\label{confstruceq}
\end{gather}

The structure equations for the isometries are obtained similarly but do not require prolongation.
The structure group is as above except $r=0$ and the structure equations are:
\begin{gather}
d\Theta_1= \Pi \wedge \Theta_2,
\qquad
d\Theta_2= \Pi \wedge \Theta_1,
\qquad
d\Theta_3= \Theta_2 \wedge \Theta_1,
\qquad
d\Pi= \kappa \Theta_2 \wedge \Theta_1.
\label{isostruceq}
\end{gather}

In both sets of structure equations the coef\/f\/icients are all constant and as a~consequence the symmetries of these
structure equations are given by the Lie group which they def\/ine (see the remarks following~\cite[Theorem 8.22]{Olver1}).
By construction the lift of our original symmetry is a~symmetry of the structure equations and must therefore be given
by the action of the Lie group which the structure equations def\/ine.

From~\eqref{isostruceq} we see that the isometries are at most a~$4$-dimensional Lie group.
By Tanaka's theory~\cite{tanak1}, the maximal dimension is reached when $\kappa=0$ and the structure equations are those of the
Heisenberg group extended by the action of a~particular strata preserving derivation of the Heisenberg algebra.
The Lie algebra of this group has the form
\begin{gather*}
[e_1, e_2]= e_3,
\qquad
[e_4, e_1]= e_2,
\qquad
[e_4, e_2]= e_1,
\end{gather*}
where $\{e_1,e_2,e_3\}$ is a~basis for the Heisenberg algebra and $e_4$ is the derivation.

The Killing form for the conformal structure equations is
\begin{gather*}
K=\left(
\begin{matrix}
0 & 0 & 0 & 0 & 0 & 0 &-7 & 0
\\
0 & 0 & 0 & 0 & 0 & 6 & 0 & 0
\\
0 & 0 & 0 & 0 & 0 & 0 & 0 & 6
\\
0 & 0 & 0 & 12 & 0 & 0 & 0 & 0
\\
0 & 0 & 0 & 0 & 4 & 0 & 0 & 0
\\
0 & 6 & 0 & 0 & 0 & 0 & 0 & 0
\\
-7 & 0 & 0 & 0 & 0 & 0 & 0 & 0
\\
0 & 0 & 6 & 0 & 0 & 0 & 0 & 0
\end{matrix}
\right).
\end{gather*}
Since $\det  K \ne 0$ $(=-3048192)$ and~$K$ is indef\/inite with signature $++++---$, the Lie algebra must be
$\mathfrak{sl}_3$, i.e., the only $8$-dimensional simple Lie algebras are $\mathfrak{sl}_3$, $\mathfrak{su}_3$ and
$\mathfrak{su}_{2,1}$, however $\mathfrak{su}_3$ and $\mathfrak{su}_{2,1}$ are ruled out by the indef\/initeness and
signature.
Alternatively one can simply compute the Lie brackets of the vector f\/ields dual to the system of one forms and check
that it isomorphic to~$\mathfrak{sl}_3$.
We should remark that Ian Anderson's MAPLE Dif\/ferentialGeometry
packages were an indispensable tool used in the
calculations outlined above.

The fact that~$\kappa$ is not present in~\eqref{confstruceq} implies that the universal cover of ${\rm SL}_2(\mathbb{R})$ and
the Heisenberg group both have ${\rm SL}_3(\mathbb{R})$ as the conformal symmetry group and consequently are conformally
equivalent, see~\cite[Proposition~2.3.2]{krugthe} and~\cite[Section~2.5]{Acap}.
We thus have the following conformal Darboux theorem.

\begin{Theorem}
All left-invariant sub-Lorentzian structures on the universal cover of ${\rm SL}_2(\mathbb{R})$ such that $\tilde h=0$ are
locally conformally equivalent to the sub-Lorentzian Heisenberg group.
\end{Theorem}

\section{Rigid example}\label{Section6}

By def\/inition, any left translation is an isometry of a~left-invariant structure and so the dimension of the isometry
group is at least~$3$, and from the previous section the dimension of the isometry group is at most~$4$.
The goal of this section is to show that the extreme case of no isometries can occur but obviously not for a~left
invariant structure.
A~relatively straight forward example comes from the geometry of second-order ODEs.

To a~given second-order ODE
\begin{gather}
u^{\prime \prime}=Q(x,u,u^{\prime}),
\label{sode}
\end{gather}
where~$Q$ is smooth, we associate three $1$-forms given by
\begin{gather}
\omega^{1}=du-pdx,
\qquad
\omega^{2}=dp-Q(x,u,p)dx,
\qquad
\omega^{3}=dx,
\label{1forms}
\end{gather}
which are regarded as one forms jet space $J^{1}(\mathbb{R}, \mathbb{R})$ with coordinates by $(x,u,p)$, where~$x$ is
the independent variable,~$u$ is the dependent variable and $p=u^{\prime}$, see~\cite{Olver1}.

In particular a~curve $\gamma (x)=(x,u(x),p(x))$ in $J^{1}(\mathbb{R}, \mathbb{R})$ def\/ines a~solution to~\eqref{sode}
if and only if $\gamma^{\ast}\omega^{i}=0$, $i=1,2$ (one can easily show that the vanishing of the two pull-backs is
equivalent to $y^{\prime}(x)=p(x)$ and in turn to $u^{\prime \prime}(x)=Q(x,u(x),u^{\prime}(x))$).
A~local dif\/feomorphism $\Phi:\mathbb{R} \times \mathbb{R\longrightarrow R}\times \mathbb{R}$ is called a~\textit{point
transformation} or \textit{point symmetry} of~\eqref{sode} if and only if it maps the graph of a~solutions
to~\eqref{sode} onto graphs of solutions to~\eqref{sode}.

Any local dif\/feomorphism $\Phi:\mathbb{R} \times \mathbb{R\longrightarrow R} \times \mathbb{R}$ can always be
\textit{prolonged} to a~local dif\/feo\-mor\-phism $\hat \Phi: J^{1}(\mathbb{R}, \mathbb{R}) \to J^{1}(\mathbb{R},
\mathbb{R})$ by setting $\hat \Phi(x,u,p)=(\tilde x, \tilde u, \tilde p)$, where
\begin{gather*}
(\tilde x, \tilde u)=\Phi(x,u)
\qquad
\text{and}
\qquad
\tilde p= \frac{d\tilde u}{d \tilde x}.
\end{gather*}
By direct calculation it follows that $\Phi$ is a~point symmetry of~\eqref{sode} if and only if there exists smooth
functions $a_{i}$, $i=1,\dots,5$, such that
\begin{gather*}
\hat{\Phi}^{\ast}\omega^{1}=a_{1}\omega^{1},
\qquad
\hat{\Phi}^{\ast}\omega^{2}= a_{2}\omega^{1}+a_{3}\omega^{2},
\qquad
\hat{\Phi}^{\ast}\omega^{3}=a_{4}\omega^{1}+a_{5}\omega^{3}.
\end{gather*}
A~classical problem in the geometric theory of ODEs is the classif\/ication of second-order ODEs with respect to point
transformations.
The fundamental result is as follows (see~\cite[Theorem~12.19]{Olver1}):

\begin{Theorem}
The point transformation symmetry group of a~second-order ordinary differential equation has dimension at most eight.
Moreover, the equation admits an eight-dimensional symmetry group if and only if it can be mapped by a~point
transformation to the linear equation $u^{\prime \prime}=0$, which has symmetry group ${\rm SL}(3)$.
\end{Theorem}

Thus the equation $u^{\prime \prime}=0$ has the maximal possible point symmetry group while at the other end of the
scale the following equation
\begin{gather}
u^{\prime \prime}=\big(\big(x+x^{2}\big)e^{u}\big)^{\prime}
\label{sodespec}
\end{gather}
has no nontrivial point symmetries (see~\cite[p.~182]{Olver1}).

Any equation~\eqref{sode} def\/ines through the forms~\eqref{1forms} a~conformal class of contact sub-Lorentzian metrics.
Indeed, let $H=\ker \omega^{1}$, $L_{1}=\ker \omega^{1}\cap \ker \omega^{2}$ and $L_{2}=\ker \omega^{1}\cap \ker \omega^{3}$.
Clearly,~$H$ is a~contact distribution that splits into the union of line bundles: $H=L_{1}\oplus L_{2}$.
Similarly as in the classical situation (cf.~\cite{Beem}) the splitting can be viewed as the f\/ield of null cones for
a~Lorentzian metric on~$H$.
Of course all such metrics are conformally equivalent.
In particular, it is seen that all solutions to~\eqref{sode} are determined by the trajectories of the null f\/ield
$\frac{\partial}{\partial x}+p\frac{\partial}{\partial u}+Q(x,u,p) \frac{\partial}{\partial p}$ spanning $L_{1}$.

\begin{Proposition}
\label{stwierdz}
Fix an equation~\eqref{sode} and let $(H,g)$ be a~sub-Lorentzian metric belonging to the conformal class of
sub-Lorentzian metrics induced by this equation.
Then any isometry of~$(H,g)$ which is isotopic to the identity is in fact a~point symmetry of the considered equation.
\end{Proposition}

\begin{proof}
Let $F:J^{1}\longrightarrow J^{1}$ be an isometry of $(H,g)$ as in the hypothesis of the proposition.
Then, obviously, $dF(H)\subset H$ and $dF(L_{i})\subset L_{i}$, $i=1,2$.
Using the above notation this last remark is equivalent to
the equations $F^{\ast}\omega^{1}=a_{1}\omega^{1}$,
$F^{\ast}\omega^{2}=a_{2}\omega^{1}+a_{3}\omega^{2}$,
$F^{\ast}\omega^{3}=a_{4}\omega^{1}+a_{5}\omega^{3}$ for some smooth functions $a_{i}$, $i=1,\dots,5$.
If (locally) we set $F=(F_{1},F_{2},F_{3})$, then the f\/irst and the third equation shows that $F_{1}$ and $F_{2}$ do not
depend on~$p$, which means that~$F$ is the f\/irst jet prolongation of a~dif\/feomorphism of the $(x,u)$-space.
This last statement is equivalent to saying that~$F$ is a~point symmetry of~\eqref{sode}.
\end{proof}

Let~$U$ be a~neighborhood of $0$ in $\mathbb{R}^{3}$ and consider any sub-Lorentzian structure $(U,H,g)$ belonging to
the conformal class of equation~\eqref{sodespec}.

\begin{Corollary}
The algebra of infinitesimal isometries of $(U,H,g)$ is trivial.
\end{Corollary}

\begin{proof}
Indeed, suppose that~$X$ is an inf\/initesimal isometry with f\/low $\psi^t$.
Since $\psi^t$ is an isometry of $(U,H,g)$ isotopic to the identity, it follows from Proposition~\ref{stwierdz} that
$\psi^t$ is a~point symmetry of~\eqref{sodespec} and therefore $\psi^t={\rm id}$.
Consequently we must have $X=0$.
\end{proof}

We remark that the construction of sub-Lorentzian structure from an ODE as above is a~particular example of a~more
general theory relating ODEs and what is sometimes called a~para-cr
structure.
On this subject we refer the reader to~\cite{HN}.

\section{Appendix}

In this appendix we would like to draw the reader's attention to some possible applications of the invariants to
non-contact cases.
Consider the simplest such case, namely the Martinet case.
Martinet sub-Lorentzian structures (of Hamiltonian type) were studied in~\cite{Gro01}.
Let $(M,H,g)$ be a~sub-Lorentzian manifold where $(H,g)$ is \textit{a Martinet sub-Lorentzian structure $($or a~metric$)$}.
That is, there exists a~hypersurface~$S$, the so-called Martinet surface, with the following pro\-per\-ties:
\begin{enumerate}\itemsep=0pt
\item[1)] $H$ is a~contact structure on $M\backslash S$,
\item[2)] $\dim (H_{q}\cap T_{q}S)=1$ for every $q\in M$,
\item[3)] the f\/ield of directions $L:S\ni q\longrightarrow L_{q}=H_{q}\cap T_{q}S$ is timelike.
\end{enumerate}

It is a~standard fact that trajectories of~$L$ are abnormal curves for the distribution~$H$.
Obviously our construction of the invariants can be carried out on the contact
sub-Lorentzian manifold $(M\backslash S,H_{|M\backslash S},g_{|M\backslash S})$.
In this way we can produce necessary conditions for two Martinet sub-Lorentzian structures to be $ts$-isometric.
More precisely, let $(M_{i},H_{i},g_{i})$ be Martinet sub-Lorentzian manifolds such that $(H_{i},g_{i})$ are
$ts$-oriented Martinet sub-Lorentzian metrics for $i=1,2$. Suppose that $\varphi:(M_{1},H_{1},g_{1})\longrightarrow
(M_{2},H_{2},g_{2})$ is a~$ts$-isometry, then since abnormal curves are preserved by dif\/feomorphisms,
$\varphi (S_{1})=S_{2}$, where $S_{i}$ is the Martinet surface for $H_{i}$, $i=1,2$. It follows that $\varphi$ induces
a~$ts$-isometry $\tilde{\varphi}=\varphi_{|M_{1}\backslash S_{1}}:(M_{1}\backslash S_{1},H_{1|M_{1}\backslash
S_{1}},g_{1|M_{1}\backslash S_{1}})\longrightarrow (M_{2}\backslash S_{2},H_{2|M_{2}\backslash
S_{2}},g_{2|M_{2}\backslash S_{2}})$.
Therefore, using results from Section~\ref{Section4} we arrive at
\begin{gather*}
\chi_{1}=\tilde{\varphi}^{\ast}\chi_{2},
\qquad
\kappa_{1}=\tilde{\varphi}^{\ast}\kappa_{2},
\qquad
\text{and}
\qquad
\tilde{h}_{1}=\tilde{\varphi}^{\ast}\tilde{h}_{2},
\end{gather*}
where $\chi_{i}$, $\kappa_{i}$, $\tilde{h}_{i}$ are the corresponding invariants for $(M_{i}\backslash
S_{i},H_{i|M_{i}\backslash S_{i}},g_{i|M_{i}\backslash S_{i}})$, $i=1,2$.

As one might expect, the invariants become singular when one approaches the Martinet surface.
Indeed, let us look at the following example.

\begin{Example}
Consider the simplest Martinet sub-Lorentzian structure, namely the f\/lat one (cf.~\cite{Gro01}).
This structure is def\/ined on $\mathbb{R}^{3}$ by the orthonormal frame
\begin{gather*}
X_{1}=\frac{\partial}{\partial x}+\frac{1}{2}y^{2}\frac{\partial}{\partial z},
\qquad
X_{2}=\frac{\partial}{\partial y}-\frac{1}{2}xy\frac{\partial}{\partial z},
\end{gather*}
where we assume $X_{1}$ (resp.\
$X_{2}$) to be a~time (resp.\
space) orientation.
The Martinet surface in this case is $S=\left\{y=0\right\}$, and we can write $H=\Span\{X_{1},X_{2}\}=\ker \omega$ for
$\omega$ def\/ined as $\omega=\frac{2}{3}\frac{1}{y}dz-\frac{1}{3}ydx+\frac{1}{3}xdy$.
Clearly, $d\omega (X_{1},X_{2})=1$, and as usual we def\/ine the Reeb f\/ield $X_{0}$ on $\mathbb{R}^{3}\backslash S$ with
equations $d\omega (X_{0},\cdot)=0$, $\omega (X_{0})=1$.
Direct computation yields
\begin{gather*}
X_{0}=-\frac{1}{y}\frac{\partial}{\partial x}+y\frac{\partial}{\partial z}.
\end{gather*}
Moreover
\begin{gather*}
[X_{2},X_{1}]=\frac{1}{y}X_{1}+X_{0},
\qquad
[X_{1},X_{0}]=0,
\qquad
[X_{2},X_{0}]=\frac{1}{y^{2}}X_{1},
\end{gather*}
from which we f\/inally obtain
\begin{gather*}
\tilde{h}=\left(
\begin{matrix}
0 & \frac{1}{2}\frac{1}{y^{2}}
\\
-\frac{1}{2}\frac{1}{y^{2}} & 0
\end{matrix}
\right),
\qquad
\chi=\frac{1}{4}\frac{1}{y^{4}},
\qquad
\text{and}
\qquad
\kappa=-\frac{5}{2}\frac{1}{y^{2}}.
\end{gather*}
\end{Example}

\subsection*{Acknowledgements}
We would like to thank the referees for their thoughtful comments and careful reading of the paper.

\pdfbookmark[1]{References}{ref}
\LastPageEnding


\begin{thebibliography}{99}
\footnotesize\itemsep=0pt

\bibitem{Agr2}
Agrachev A.A., Barilari D., Sub-{R}iemannian structures on 3{D} {L}ie groups,
  \href{http://dx.doi.org/10.1007/s10883-012-9133-8}{\textit{J.~Dyn. Control Syst.}} \textbf{18} (2012), 21--44,
  \href{http://arxiv.org/abs/#2}{arXiv:1007.4970}.

\bibitem{Agr1}
Agrachev A.A., El~Alaoui El-H.Ch., Gauthier J.-P., Sub-{R}iemannian metrics on
  {${\mathbb R}^3$}, in Geometric Control and Non-Holonomic Mechanics ({M}exico
  {C}ity, 1996), \textit{CMS Conf. Proc.}, Vol.~25, Amer. Math. Soc.,
  Providence, RI, 1998, 29--78.

\bibitem{Beem}
Beem J.K., Ehrlich P.E., Easley K.L., Global {L}orentzian geometry,
  \textit{Monographs and Textbooks in Pure and Applied Mathematics}, Vol.~202,
  2nd ed., Marcel Dekker, Inc., New York, 1996.

\bibitem{BerGic}
Berestovskii V.N., Gichev V.M., Metrized semigroups, \href{http://dx.doi.org/10.1023/B:JOTH.0000008737.57612.5d}{\textit{J.~Math. Sci.}}
  \textbf{19} (2004), 10--29.

\bibitem{Boar}
Boarotto F., Conformal equivalence of sub-Riemannian 3D contact structures on
  Lie groups, \href{http://arxiv.org/abs/1412.2358}{arXiv:1412.2358}.

\bibitem{Acap}
{\v{C}}ap A., Automorphism groups of parabolic geometries, \textit{Rend. Circ.
  Mat. Palermo}  (2005), suppl., 233--239.

\bibitem{Agr0}
El~Alaoui El-H.Ch., Gauthier J.-P., Kupka I., Small sub-Riemannian balls in
  $\mathbb{R}^{3}$, \href{http://dx.doi.org/10.1007/BF02269424}{\textit{J.~Dyn. Control Syst.}} \textbf{2} (1996), 359--421.

\bibitem{FalGor}
Falbel E., Gorodski C., Sub-{R}iemannian homogeneous spaces in dimensions {$3$}
  and {$4$}, \href{http://dx.doi.org/10.1007/BF00181566}{\textit{Geom. Dedicata}} \textbf{62} (1996), 227--252.

\bibitem{g2}
Grochowski M., Normal forms of germs of contact sub-{L}orentzian structures on
  {${\mathbb R}^3$}. {D}if\/ferentiability of the sub-{L}orentzian distance
  function, \href{http://dx.doi.org/10.1023/A:1025696302287}{\textit{J.~Dynam. Control Systems}} \textbf{9} (2003), 531--547.

\bibitem{bcc2004}
Grochowski M., On the {H}eisenberg sub-{L}orentzian metric on {${\mathbb
  R}^3$}, in Geometric Singularity Theory, \href{http://dx.doi.org/10.4064/bc65-0-4}{\textit{Banach Center Publ.}},
  Vol.~65, Polish Acad. Sci., Warsaw, 2004, 57--65.

\bibitem{g7}
Grochowski M., Reachable sets for the {H}eisenberg sub-{L}orentzian structure
  on {${\mathbb R}^3$}. {A}n estimate for the distance function,
  \href{http://dx.doi.org/10.1007/s10450-006-0378-y}{\textit{J.~Dyn. Control Syst.}} \textbf{12} (2006), 145--160.

\bibitem{g6}
Grochowski M., Properties of reachable sets in the sub-{L}orentzian geometry,
  \href{http://dx.doi.org/10.1016/j.geomphys.2009.03.011}{\textit{J.~Geom. Phys.}} \textbf{59} (2009), 885--900.

\bibitem{Gro01}
Grochowski M., Normal forms and reachable sets for analytic {M}artinet
  sub-{L}orentzian structures of {H}amiltonian type, \href{http://dx.doi.org/10.1007/s10883-011-9110-7}{\textit{J.~Dyn. Control
  Syst.}} \textbf{17} (2011), 49--75.

\bibitem{Gro02}
Grochowski M., Reachable sets for contact sub-{L}orentzian structures on
  {${\mathbb R}^3$}. {A}pplication to control af\/f\/ine systems on {${\mathbb
  R}^3$} with a scalar input, \href{http://dx.doi.org/10.1007/s10958-011-0464-x}{\textit{J.~Math. Sci.}} \textbf{177} (2011),
  383--394.

\bibitem{groGlob}
Grochowski M., Remarks on global sub-{L}orentzian geometry, \href{http://dx.doi.org/10.1007/s13324-013-0059-4}{\textit{Anal. Math.
  Phys.}} \textbf{3} (2013), 295--309.

\bibitem{GMW}
Grochowski M., Medvedev A., Warhurst B., Classif\/ication of 3-dimensional
  contact left-invariant sub-Lorentzian structures, {i}n preparation.

\bibitem{Grong}
Grong E., Vasil'ev A., Sub-{R}iemannian and sub-{L}orentzian geometry on {${\rm
  SU}(1,1)$} and on its universal cover, \href{http://dx.doi.org/10.3934/jgm.2011.3.225}{\textit{J.~Geom. Mech.}} \textbf{3}
  (2011), 225--260, \href{http://arxiv.org/abs/0910.0945}{arXiv:0910.0945}.

\bibitem{HN}
Hill C.D., Nurowski P., Dif\/ferential equations and para-{CR} structures,
  \textit{Boll. Unione Mat. Ital.} \textbf{3} (2010), 25--91,
  \href{http://arxiv.org/abs/0909.2458}{arXiv:0909.2458}.

\bibitem{Huang}
Huang T., Yang X., Geodesics in the {H}eisenberg group {$H^n$} with a
  {L}orentzian metric, \href{http://dx.doi.org/10.1007/s10883-012-9156-1}{\textit{J.~Dyn. Control Syst.}} \textbf{18} (2012),
  479--498.

\bibitem{markin}
Korolko A., Markina I., Nonholonomic {L}orentzian geometry on some {${\mathbb
  H}$}-type groups, \href{http://dx.doi.org/10.1007/s12220-009-9088-5}{\textit{J.~Geom. Anal.}} \textbf{19} (2009), 864--889,
  \href{http://arxiv.org/abs/0809.4450}{arXiv:0809.4450}.

\bibitem{markin2}
Korolko A., Markina I., Geodesics on {${\mathbb H}$}-type quaternion groups
  with sub-{L}orentzian metric and their physi\-cal interpretation,
  \href{http://dx.doi.org/10.1007/s11785-010-0079-0}{\textit{Complex Anal. Oper. Theory}} \textbf{4} (2010), 589--618,
  \href{http://arxiv.org/abs/1004.1508}{arXiv:1004.1508}.

\bibitem{krugthe}
Kruglikov B., The D., The gap phenomenon in parabolic geometries,
  \href{http://arxiv.org/abs/1303.1307}{arXiv:1303.1307}.

\bibitem{LiuSuss}
Liu W., Sussman H.J., Shortest paths for sub-{R}iemannian metrics on rank-two
  distributions, \textit{Mem. Amer. Math. Soc.} \textbf{118} (1995),
  x+104~pages.

\bibitem{Olver1}
Olver P.J., Equivalence, invariants, and symmetry, \href{http://dx.doi.org/10.1017/CBO9780511609565}{Cambridge University Press},
  Cambridge, 1995.

\bibitem{tanak1}
Tanaka N., On dif\/ferential systems, graded {L}ie algebras and pseudogroups,
  \textit{J.~Math. Kyoto Univ.} \textbf{10} (1970), 1--82.

\end{thebibliography}
\end{document}